\documentclass[11pt]{article}
\usepackage[T1]{fontenc}
\usepackage[utf8]{inputenc}
\usepackage{amsmath}
\usepackage{amsfonts, amssymb}
\usepackage[dvipsnames]{xcolor}
\usepackage{graphicx}
\newtheorem{theo}{Theorem}[section]
\newtheorem{lem}[theo]{Lemma}
\newtheorem{cor}[theo]{Corollary}
\newtheorem{prop}[theo]{Proposition}

\newtheorem{conj}[theo]{Conjecture}

\newcommand{\mysection}[1]{\section{#1} \setcounter{equation}{0}}
\newcommand{\proof}{{\sc Proof.} \quad}
\newcommand{\proofc}{{\sc Proof} \ }
\newcommand{\be}{\begin{equation} \label}
\newcommand{\ee}{\end{equation}}
\newcommand{\bea}{\begin{eqnarray}\label}
\newcommand{\eea}{\end{eqnarray}}
\newcommand{\bas}{\begin{eqnarray*}}
\newcommand{\eas}{\end{eqnarray*}}
\newcommand{\bit}{\begin{itemize}}
\newcommand{\eit}{\end{itemize}}
\newcommand{\qed}{\hfill$\Box$ \vskip.2cm}
\newcommand{\nn}{\nonumber}
\newcommand{\R}{\mathbb{R}}
\newcommand{\N}{\mathbb{N}}
\newcommand{\pO}{\partial\Omega}
\newcommand{\bom}{\overline{\Omega}}
\newcommand{\eps}{\varepsilon}

\newcommand{\io}{\int_\Omega}

\newcommand{\Mint}{- \hspace*{-4mm} \int}

\newcommand{\abs}{\\[5pt]}
\newcommand{\tm}{T_{max}}

\newcommand{\Mstar}{M_\star}
\newcommand{\MStar}{M^\star}
\newcommand{\mstar}{m_\star}
\newcommand{\mStar}{m^\star}
\newcommand{\uw}{\underline{w}}
\newcommand{\ow}{\overline{w}}
\newcommand{\om}{\overline{m}}
\newcommand{\oy}{\overline{y}}
\newcommand{\win}{w_{in}}
\newcommand{\wout}{w_{out}}
\newcommand{\F}{{\mathcal{F}}}
\newcommand{\D}{{\mathcal{D}}}
\newcommand{\C}{{\mathcal{C}}}
\newcommand{\Om}{\Omega}
\newcommand{\dOm}{\partial\Omega}
\newcommand{\Ombar}{\bom}
\newcommand{\norm}[2][]{\|#2\|_{#1}}
\newcommand{\normm}[2]{\|#2\|_{#1}}
\newcommand{\Lom}[1]{L^{#1}(\Om)}
\newcommand{\na}{\nabla}

\newcommand{\f}[2]{\frac{#1}{#2}}
\newcommand{\set}[1]{\left\{#1\right\}}
\newcommand{\ww}{\underline{w}}
\newcommand{\kl}[1]{\left(#1\right)}
\usepackage{newunicodechar}
\newunicodechar{Δ}{\Delta}
\newunicodechar{∇}{\nabla}
\newunicodechar{ν}{\nu}
\newunicodechar{φ}{\varphi}
\newunicodechar{Φ}{\Phi}
\newunicodechar{∞}{\infty}
\newunicodechar{ψ}{\psi}
\newunicodechar{τ}{\tau}
\newunicodechar{α}{\alpha}
\newunicodechar{π}{\pi}
\newunicodechar{ℝ}{\mathbb{R}}
\begin{document}
\enlargethispage{20mm}
\title{A double critical mass phenomenon in a\\
no-flux-Dirichlet Keller-Segel system}
\author{
Jan Fuhrmann\footnote{j.fuhrmann@fz-juelich.de}\\
{\small J\"ulich Supercomputing Centre, Forschungszentrum J\"ulich}\\
{\small 52428 J\"ulich, Germany}\\
{\small and Frankfurt Institute for Advanced Studies}\\
{\small 60438 Frankfurt/Main, Germany}
\and
Johannes Lankeit\footnote{lankeit@ifam.uni-hannover.de}\\
{\small Leibniz Universität Hannover, Institut f\"ur Angewandte Mathematik,}\\
{\small Welfengarten 1, 30167 Hannover, Germany} 
\and
Michael Winkler\footnote{michael.winkler@math.uni-paderborn.de}\\
{\small Institut f\"ur Mathematik, Universit\"at Paderborn,}\\
{\small 33098 Paderborn, Germany} 
}
\date{}
\maketitle

\vspace*{-6mm}
\begin{abstract}
\noindent 
Derived from a biophysical model for the motion of a crawling cell, the evolution system
\bas
	\left\{ \begin{array}{l}
	u_t=\Delta u -\nabla \cdot (u\nabla v), \\[1mm]
	0 = \Delta v -kv+u,
	\end{array} \right.
	\qquad  \qquad (\star)
\eas
is investigated in a finite domain $\Omega\subset\R^n$, $n\ge 2$, with $k\ge 0$.
Whereas a comprehensive literature is available for cases in which ($\star$) describes chemotaxis-driven population dynamics
and hence is accompanied by homogeneous Neumann-type boundary conditions for both components,
the presently considered modeling context, besides yet requiring the flux $\partial_\nu u - u\partial_\nu n$ to vanish on $\pO$,
inherently involves homogeneous Dirichlet boundary conditions for the attractant $v$, 
which in the current setting corresponds to the cell's cytoskeleton being free of pressure at the boundary.\abs
This modification in the boundary setting is shown to go along with a substantial change with respect to the potential to 
support the emergence of singular structures:
It is, inter alia, revealed that in contexts of radial solutions in balls there exist {\em two} critical mass levels,
distinct from each other whenever $k>0$ or $n\ge 3$,
that separate ranges within which $(i)$ all solutions are global in time and remain bounded, 
$(ii)$ both global bounded and exploding solutions exist, or $(iii)$ all nontrivial solutions blow up in finite time.
While critical mass phenomena distinguishing between regimes of type $(i)$ and $(ii)$ belong to the well-understood
characteristics of ($\star$) when posed under classical no-flux boundary conditions in planar domains, 
the discovery of a distinct secondary critical mass level related to the occurrence of $(iii)$ seems to have no nearby precedent.\abs
In the planar case with the domain being a disk, the analytical results are supplemented with some numerical illustrations,
and it is discussed how the findings can be interpreted biophysically for the situation of a cell on a flat substrate.\abs
\end{abstract}

%

%
%
%
%
%
%
%
%
\section{Introduction}\label{intro}
{\bf A Keller-Segel type model for crawling keratocytes.} \quad
This study is concerned with the cross-diffusion problem
\be{01}
	\left\{ \begin{array}{ll}
	u_t=\Delta u -\nabla \cdot (u\nabla v),
	\qquad & x\in\Omega, \ t>0, \\[1mm]
	0 = \Delta v -kv+u,
	\qquad & x\in\Omega, \ t>0, \\[1mm]
	u(x,0)=u_0(x),
	\qquad & x\in\Omega,
	\end{array} \right.
\ee
in a bounded domain $\Omega\subset\R^n$, $n\ge 2$.
During the past decades, this system has received noticeable interest when used as a parabolic-elliptic simplification of the
celebrated Keller-Segel model to describe collective behavior in microbial populations with movement chemotactically 
biased by a chemical signal, and hence typically found accompanied by no-flux boundary conditions in the literature
(\cite{KS}, \cite{horstmann03}, \cite{LWsurvey}).\abs
In contrast to this, the context to be considered in the present paper necessitates to supplement (\ref{01}) by the requirements
\be{02}
	\frac{\partial u}{\partial\nu}- u\frac{\partial v}{\partial\nu}=v=0,
	\qquad x\in\pO,
\ee
on the boundary of the domain $\Omega\subset\R^n$, as intrinsically linked to the role which, quite independently of the 
above, (\ref{01}) plays when
derived from a biomechanical model for a single crawling keratocyte, or rather a keratocyte fragment, that has been
introduced in \cite{BLATM} for space dimension $n=2$. These fragments are similar to lamellipodia, i.e., very flat structures, 
and can in good approxomation be described as two-dimensional entities. The computational model presented in \cite{BLATM} was reduced and analyzed in \cite{BFR}, and similar models in one space dimension have been investigated in, e.g., \cite{RPT}. 
From the physical model in \cite{BLATM}, a reduced free boundary problem has been derived in \cite{BFR} by combining bulk and 
shear components of the stress in the actin gel in a phenomenological way, allowing for the stress tensor to be represented 
as a scalar multiple of the identity matrix. This step used the fact that cytoskeleton gels are rather unusual viscoelastic fluids 
with the stress not being shear dominated. This led to a free boundary problem for two variables, in our context named $v$ 
for the stress in the cytoskeleton and $u$ for the density of myosin motor proteins. The latter actively 
generate stress by binding to and pulling on the actin filaments 
constituting the cytoskeleton meshwork.\abs 
The first equation in \eqref{01} is thus interpreted as a diffusion-advection equation for the concentration of myosin molecules which
are either freely diffusing inside the cytoplasm or are bound to the actin gel and hence convected with the velocity $\nabla v$ which
is the divergence of the stress tensor, $v \mathbb{I}$. The second equation describes the force balance in the actin gel with the term
$u$ representing the actively generated stress due to the myosin motors, which is assumed to be proportional to the density of these motors. The term $-kv$ models the dissipation of stress via traction with the substrate to which the actin gel is linked by adhesion
molecules. The distribution of these adhesions is supposed to be uniform and constant in time for a resting cell. Moreover, the second
equation being elliptic assumes that stresses equilibrate on a much faster time scale than the motion of the actin gel, indicated by
very low Deborah numbers reported for moving, let alone resting cells \cite{RFJVM}. This simply means that the gel behaves more 
like a viscous fluid than an elastic solid on the relevant time scale. The parameter $k$ is the typical stress stored in the actin 
gel relative to the typical stress generated by myosin motors. The second parameter present in the model is the size of the 
domain $\Omega$ which is measured in multiples of $\sqrt{k}\mathcal{L}$, where $\mathcal{L}$ is the viscous length of the actin 
gel which describes how far the locally generated stress acts through the network before being dissipated away. 
It is defined as square root of the ratio of the viscosity and the traction coefficient.\abs
Whereas both \cite{BLATM} and \cite{BFR} were interested in traveling wave solutions to their respective free boundary problems to
describe steady cell motion, we will focus here on the behavior of steady states and the possibility of finite time blow up. Steady
state solutions clearly correspond to a resting cell although we should mention that stationarity in \eqref{01} does not imply that
there is no motion inside the cell; recall that the velocity of the actin gel is $\nabla v$. More strikingly, solutions blowing up 
in finite time are interpreted as the cell being physically disrupted by too much contractile activity of myosin motors as
represented by a large total myosin mass $m = \int_\Omega u$ which is obviously a conserved quantity for (\ref{01})-(\ref{02}). 
While the bifurcation from rest to motion at subcritical values of $m$ described in \cite{BFR} refers to a dynamic instability 
of the free boundary problem modeling a potentially motile cell switching from rest to directed motion, blow up of solutions 
for large $m$ in system (\ref{01}) with fixed boundary relates to the observed disruption of immobile cells upon variations of myosin
activity or adhesion strength as has been seen experimentally (\cite{BWHI}, cf. e.g. \cite{Schr2019} for mechanism of fragmentation 
of actin filaments by myosin generated forces). Mechanical breakdown due to enhanced myosin activity and concomitant concentration of
myosin is also associated with physiological processes such as programmed cell death, or apoptosis, as described in \cite{Croft}.\abs
To rule out possible issues of self intersection of the moving boundary as mechanism for the break down of solutions we fixed the shape
of the domain $\Omega$ occupied by the cell. Physically, this may be achieved by letting the cell sit on a particularly sticky
substrate or by providing it with an adhesive patch of substrate of a given shape $\Omega$ and making the surrounding region, 
viz. $\mathbb{R}^2\backslash\Omega$, particularly hostile by coating with adverse substances or no coating at all. Keeping the stress-free boundary condition $v=0$ and the no-flux condition for the myosin molecules from the original model (\cite{BLATM}), we finally arrive at (\ref{01})-(\ref{02}) which differs from the classical parabolic-elliptic Keller-Segel system most significantly in the boundary conditions. The peculiar condition $v=0$ on $\pO$ arises from the fact that myosin motors at the boundary are not supposed to generate stress since there is nothing outside the cell to be pulled against. There is no contradiction in the cytoskeleton gel's velocity being different from zero at the boundary. In fact, in a resting cell, actin is polymerized at the boundary, leading on averaege to a radial expansion of the cytoskeleton, which is counteracted by the actin gel constantly moving toward the center where the actin filaments are depolymerized. This retrograde flow means that the gel moves away from the boundary at non-zero velocity.\abs
{\bf Detecting explosion-related dichotomies in Keller-Segel systems.} \quad
Over the past decades, significant effort in the analysis of chemotaxis problems has been directed towards
excluding (e.g. \cite{NSY}) or detecting blow-up (\cite{JL,HV,nagai2001}) and the study of additional qualitative properties (e.g. \cite{senba_suzuki,mizoguchi_souplet,velazquez,blanchet_carrillo_masmoudi,blanchet_dolbeault_perthame}) in \eqref{01} and related variants, e.g.~further simplified like in \cite{JL}, or rather fully parabolic and hence more complex. Among the apparently most striking characteristics of such Keller-Segel systems, the literature has identified situations in which the occurrence of blow-up depends on the size of the conserved total mass
$\io u$ in a crucial manner. Specifically, when posed along with homogeneous Neumann boundary conditions for both components in planar bounded domains $\Omega$, \eqref{01} with arbitrary $k>0$ is known to exhibit a sharp and well-understood critical mass phenomenon in the sense that whenever $0\le u_0$ is sufficiently regular with $\io u_0 < 4\pi$, an associated initial-boundary value problem with $u|_{t=0}=u_0$ admits a globally defined bounded solution, whereas for any $m>4\pi$ one can find smooth initial data with $\io u_0=m$ such that the corresponding solution blows up in finite time (\cite{nagai2001}); a restriction to radially symmetric solutions in balls
increases this separating mass level to the value $8\pi$ (\cite{nagai2001}). Similar dichotomies have been detected in Neumann problems for further parabolic-elliptic and for fully parabolic relatives of \eqref{01} 
\cite{nagai1995,biler_nadzieja_I,HV,NSY}; cf.~also 
\cite{blanchet_dolbeault_perthame, wei_dongyi} for some related findings for Cauchy problems on the whole plane $\Omega=\R^2$).\abs
{\bf A secondary critical mass phenomenon enforced by Dirichlet conditions for $v$. Main results.} \quad
The present study will now reveal that when considered along with the boundary conditions in \eqref{02}, the system \eqref{01} may gain a further dynamical facet that is linked to the presence of a secondary, and apparently yet undiscovered, critical mass phenomenon.\abs
To appropriately formulate and embed our findings in this regard, let us first summarize some fundamental properties thereof, as can readily be verified upon straightforward adaptation of arguments known from the literature (cf.~e.g.~\cite{suzuki2013} for Part i), \cite{suzuki2013}, \cite{nagai2001} for Part ii), and \cite{nagai2000} for Part iii)):\abs
{\bf Theorem A} \quad
{\it
  Let $n\ge 2$ and $\Omega\subset\R^n$ be a bounded domain with smooth boundary, and let $k\ge 0$.\abs
  i) \ If $n=2$ and $u_0\in C^0(\bom)$ is nonnegative with
  \bas
	\io u_0  < 8\pi,
  \eas
  then \eqref{01}-\eqref{02} 
  possesses a global classical solution $(u,v)$ which is bounded in the sense that there exists $C>0$ such that
  \be{A1}
	\|u(\cdot,t)\|_{L^\infty(\Omega)} \le C
	\qquad \mbox{for all } t>0.
  \ee
  ii) \ If $n=2$, then for all $m>8\pi$ there exists some nonnegative $u_0\in C^0(\bom)$ with $\io u_0=m$ 
  such that the corresponding solution
  of (\ref{01})-(\ref{02}) blows up in finite time in the sense specified in Proposition \ref{prop33} below.
  Here, if $\Omega=B_R(0)$ with some $R>0$, then $u_0$ can be chosen to be radially symmetric with respect to $x=0$.\abs
  iii) \ In the case $n\ge 3$ and if $\Omega$ is star-shaped, for all $m>0$ one can find nonnegative $u_0\in C^0(\bom)$ 
  with $\io u_0=m$, radially symmetric if 
  $\Omega$ is a ball, such that the solution of (\ref{01})-(\ref{02}) blows up.\abs
}
As a direct consequence for the general, not necessarily radial case, this implies the following essentially well-known statement identifying the number $8\pi$ as a $k$-independent critical mass in (\ref{01})-(\ref{02}) when $n=2$, whereas if $n\ge 3$ then a corresponding critical mass phenomenon seems absent:\abs
{\bf Corollary B} \quad
{\it
  Let $n\ge 2$, $\Omega\subset\R^n$ be a bounded domain with smooth boundary, and $k\ge 0$.
  Then 
  \bea{Ms}
	\Mstar(\Omega,k):=\inf \bigg\{
	m>0 &\bigg|& \mbox{There exists some nonnegative } u_0\in C^0(\bom) \mbox{ with } \io u_0=m \nn\\
	& & \mbox{ such that the solution of \eqref{01}-\eqref{02} blows up} \bigg\}
  \eea
  is well-defined and satisfies
  \be{B.1}
	 \Mstar(\Omega,k)= \left\{ \begin{array}{ll}
	8\pi \qquad & \mbox{if } n=2, \\[1mm]
	0 & \mbox{if } n\ge 3.
	\end{array} \right.
  \ee
}

\vspace*{2mm}
Now the first of our main results identifies a secondary mass threshold which, as can already be stated at this stage,
at least in the case $n\ge 3$ indeed differs from the value $\Mstar(\Omega,k)=0$.
\begin{theo}\label{theo41}
  Let $n\ge 2$ and $\Omega\subset\R^n$ be a bounded domain with smooth boundary which is strictly star-shaped with respect
  to $0\in\Omega$ in the sense that
  \be{41.1}
	\gamma:=\inf_{x\in\partial\Omega} x\cdot \nu(x) >0.
  \ee
  Then for all $k\ge 0$,
  \bea{MS}
	\MStar(\Omega,k):=\inf \bigg\{
	m>0 &\bigg|&  \mbox{For all nonnegative } u_0\in C^0(\bom) \mbox{ with } \io u_0=m, \nn\\
	& & \mbox{ the solution of \eqref{01}-\eqref{02} blows up} \bigg\}
  \eea
  is well-defined and finite with
  \be{41.2}
	8\pi \le \MStar(\Omega,k) \le \frac{4|\pO|}{\gamma} + 2k|\Omega|
	\qquad \mbox{if } n=2
  \ee
  and
  \be{41.3}
	0< \MStar(\Omega,k) \le \frac{2n|\pO|}{\gamma} + 2k|\Omega|
	\qquad \mbox{if } n\ge 3.
  \ee
\end{theo}
In two-dimensional domains, however, the situation will turn out to be more subtle, involving a crucial qualitative 
dependence on whether or not the parameter $k$ is positive.
As a first step toward revealing this, let us concentrate on the 
special situation when $\Omega$ is a ball, in which the above enables us to rather explicitly estimate this secondary
critical mass, and to thereby detect, in particular, coincidence of both mass thresholds in the planar case 
when $k=0$ in such geometries.
\begin{cor}\label{cor411}
  Let $n\ge 2$, $R>0$ and $\Omega=B_R(0)\subset\R^n$. Then for all $k\ge 0$,
  \bas
	8\pi =\Mstar(B_R(0),k) \le \MStar(B_R(0),k) \le 8\pi + 2k\pi R^2
	\qquad \mbox{if } n=2
  \eas
  and
  \bas
	0= \Mstar(B_R(0),k) < \MStar(B_R(0),k) \le \frac{2\omega_n R^n}{n} + 2nk \omega_n R^{n-2}
	\qquad \mbox{if } n\ge 3,
  \eas
  where $\omega_n$ denotes the $(n-1)$-dimensional measure of the unit sphere $\partial B_1(0)$.
  In particular, for $k=0$,
  \bas
	\Mstar(B_R(0),0)=\MStar(B_R(0),0)=8\pi
	\quad \mbox{for all } R>0
	\qquad \mbox{if } n=2.
  \eas
\end{cor}
On further specializing the setup by resorting henceforth to radially symmetric solutions
in balls $\Omega=B_R(0)\subset \R^n$, $n\ge 2$, $R>0$, emanating from initial data in the space
$C^0_{rad}(\bom):=\{\varphi\in C^0(\bom) \ | \ \varphi \mbox{ is radially symmetric with respect to } x=0\}$,
we can rephrase part of Theorem A as follows.\abs
{\bf Corollary C} \quad
{\it
  Let $n\ge 2, R>0$, and $\Omega=B_R(0)\subset\R^n$, and let $k\ge 0$. Then
  \bea{ms}
	\mstar(\Omega,k)
	:=\inf \bigg\{
	m>0 &\bigg|&  \mbox{There exists some nonnegative } u_0\in C^0_{rad}(\bom) \mbox{ with } \io u_0=m \nn\\[0mm]
	& & \mbox{ such that the solution of (\ref{01})-(\ref{02}) blows up} \bigg\}
  \eea
  is well-defined with
  \bas
	\mstar(n,R,k)=\Mstar(B_R(0),k)=\left\{ \begin{array}{ll}
	8\pi \qquad & \mbox{if } n=2, \\[1mm]
	0 & \mbox{if } n\ge 3.
	\end{array} \right.
  \eas
}

Now the second of our main results makes sure that a corresponding secondary mass threshold, defined in the spirit
of Theorem \ref{theo41}, plays the role of a genuinely new critical mass for radial solutions not only when $n\ge 3$
and $k\ge 0$, but also when $n=2$ and $k>0$ is arbitrary, thus complementing the outcome of Corollary \ref{cor411}
in quite a sharp manner:
\begin{theo}\label{theo43}
  Let $n\ge 2$, $R>0$, and $\Omega=B_R(0)\subset\R^n$.
  Then for all $k\ge 0$,
  \bea{mS}
	\mStar(\Omega,k):=\inf \bigg\{
	m>0 &\bigg|&  \mbox{For all nonnegative } u_0\in C^0_{rad}(\bom) \mbox{ with } \io u_0=m, \nn\\
	& & \mbox{ the solution of (\ref{01})-(\ref{02}) blows up} \bigg\}
  \eea
  satisfies
  \be{43.1}
	\Mstar(B_R(0),k)=\mstar(n,R,k) \le \mStar(n,R,k) \le \MStar(B_R(0),k).
  \ee
  Moreover,
  \be{43.2}
	\mstar(2,R,0)=\mStar(2,R,0)=8\pi,
  \ee
  but
  \be{43.3}
	8\pi=\mstar(2,R,k) < \mStar(2,R,k)
	\qquad \mbox{for all } k>0,
  \ee
  and apart from that,
  \be{43.4}
	0=\mstar(n,R,k)<\mStar(n,R,k) 
	\quad \mbox{for all } k\ge 0
	\qquad \mbox{if } n\ge 3.
  \ee
\end{theo}

For the special case $k=0$, the finiteness of $\MStar$ 	(in $n$-dimensional balls, $n\ge 2$, but for possibly nonradial $u_0$) 
was already observed in \cite{biler_III} and that of $\mStar$ in 
\cite{biler_hilhorst_nadzieja_II}. It is remarkable that the values of $\mstar$ and $\mStar$, which coincide for $k=0$ and $n=2$, differ for positive $k$. In this sense linear signal degradation affects the blow-up affinity of \eqref{01} and makes it possible to find 
two separate critical masses in the same system. 
\mysection{Local existence and extensibility}
Let us first adapt an essentially well-established contraction-based reasoning to see that similar to its no-flux type relative, the problem \eqref{01}-\eqref{02} admits local smooth solutions which can cease to exist within finite time only when becoming unbounded with respect to the $L^\infty$ norm in their first component.
\begin{prop}\label{prop33}
  Let $n\ge 2$ and $\Omega\subset \R^n$ be a bounded domain with smooth boundary, let $k\ge 0$, and suppose that $u_0\in C^0(\bom)$ is nonnegative. 
  Then there exist $\tm\in (0,\infty]$ and a uniquely determined pair $(u,v)$ of nonnegative functions
  \be{reg}
	\left\{ \begin{array}{l}
	u\in C^0(\bom\times [0,\tm)) \cap C^{2,1}(\bom\times (0,\tm))
	\qquad \mbox{and} \\[1mm]
	v\in C^{2,0}(\bom\times (0,\tm))
	\end{array} \right.
  \ee
  which solve (\ref{01})-(\ref{02}) classically in $\Omega\times (0,\tm)$, and which are such that
  \be{ext}
	\mbox{if $\tm<\infty$, then $(u,v)$ blows up at $t=\tm$},
  \ee
  where we say that $(u,v)$ blows up at $t=\tm$ if and only if
  $\limsup_{t\nearrow\tm} \|u(\cdot,t)\|_{L^\infty(\Omega)}=\infty$.\abs
  Furthermore,
  \be{mass}
	\|u(\cdot,t)\|_{L^1(\Omega)} = \io u_0
	\qquad \mbox{for all } t\in (0,\tm).
  \ee
\end{prop}
\newcommand{\ubar}{\overline{u}}
\proof
We fix some $p>n$ and let $M:=\norm[\Lom p]{u_0}+1$. 
With $T>0$ to be determined later, we set 
\[
 X_{M,T} := \set{u\in C^0([0,T];\Lom p) \mid \norm[L^\infty((0,T);\Lom p)]{u}\le M,\quad u(\cdot,0)=u_0}.
\]
Given any $\ubar\in X_T:=C^0([0,T];\Lom p)$, for $t\in (0,T)$ letting $v(\cdot,t) \in 
W_0^{1,2}(\Omega)$ 
denote the weak solution of the Dirichlet problem for
$0=Δv(\cdot,t)-kv(\cdot,t)+\ubar(\cdot,t)$ we obtain a function $v=v(\ubar)\in 
C^0([0,T];W^{2,p}(\Om)\cap W_0^{1,p}(\Om))$ and note that due to our choice of $p$, elliptic 
regularity theory (see e.g. \cite[Thm. 37,I]{Miranda}) and a Sobolev embedding, we can find 
$c_1>0$ such that 
\[
 \normm{C^0([0,T];\Lom{∞})}{\na v(\ubar)} \le c_1 \normm{C^0([0,T];\Lom p)}{\ubar} \qquad \text{ 
for all } \ubar\in X_T.
\]
According to \cite[Thm. VI.39]{Lieberman}, for each $v(\ubar)$, $\ubar\in X_{M,T}$, the problem 
\[
 u_t=\na\cdot(\na u-u\na v(\ubar))\;\text{ in } \Om\times(0,T),\quad (\na u-u\na v(\ubar))\cdot ν=0\;\text{ on } \dOm\times(0,T),\quad u(\cdot,0)=u_0\;\text{ in } \Om,
\]
has a unique solution $u\in V_2=\set{u\in L^\infty((0,T);L^2(\Om))\mid \na u \in 
L^2(\Om\times(0,T))}$ which is nonnegative and bounded by some $c_2(M)$ in $\Om\times[0,T]$ 
(\cite[Thm. VI.40]{Lieberman}) and Hölder-continuous in $\Ombar\times(0,T)$ (\cite[Thm. 1.3 and Remark 1.3]{porzio-vespri}). We denote this solution by $Φ(\ubar)$, thus defining a mapping $Φ\colon X_{M,T}\to 
X_T$. 
For arbitrary $t\in (0,T)$, $h_1\in(0,T-t)$, $h_2\in(0,T-t-h_1)$, we let $\psi\equiv 1$ 
on $[0,t)$, $\psi\equiv 0$ on $(t+h_1,T)$ and linearly interpolated between $t$ and $t+h_1$.
Given $u_1, u_2\in X_{M,T}$, we then let
\[
\varphi(x,\tau):= \frac1{h_2} \int_{\tau}^{\tau+h_2} (\Phi(u_1)-\Phi(u_2))^{p-1} (x,s) ds 
\cdot \psi(\tau), \qquad x\in \Om, τ\in(0,T),
\]
and use this regularized version of $(\Phi(u_1)-\Phi(u_2))^{p-1}$ as test function in the 
difference of the definitions of weak solutions (cf.  \cite[p. 136]{Lieberman}) for $Φ(u_1)$ and 
$Φ(u_2)$. After successively taking $h_1\to 0$ and $h_2\to 0$ and 
several applications of Young's inequality we find that with some $c_3>0$, 
\[
 \f1p \io ((Φ(u_1)-Φ(u_2)) (t))^p \le c_3 (1+M^p) \int_0^t\io (Φ(u_1)-Φ(u_2))^p + c_3 c_2^p(M)\int_0^{t} \io |\na (v_1-v_2)|^p
\]
holds for every $t\in(0,T)$, $u_1,u_2\in X_{M,T}$. 
Therefore, by a Grönwall-type argument we find that with some $c_4>0$,
\[
 \norm[\Lom p]{Φ(u_1)(t)-Φ(u_2)(t)}^{p} \le c_4(e^{c_4t}-1)\norm[L^\infty((0,T);\Lom p]{\na v_1-\na v_2}^p \le c_1 c_4(e^{c_4T}-1)\norm[L^\infty((0,T);\Lom p]{u_1-u_2}^{p}
\]
is satisfied for all $u_1, u_2\in X_{M,T}$ and all $t\in(0,T)$. Upon suitably small choice of $T$, the map $Φ\colon X_{M,T}\to X_{M,T}$ becomes a contraction. Banach's theorem hence entails the existence of a fixed point $u=Φ(u)$, unique within $X_{M,T}$, whose further regularity follows from successive applications of \cite[Thm. 6.6]{GT}, \cite[Thm 1.1]{lieberman_paper} and \cite[Thm. IV.5.3]{LSU}. 
The extensibility criterion \eqref{ext} is a consequence of the exclusive dependence of $T$ on $M$, and hence 
on $\norm[\Lom \infty]{u_0}$, whereas (\ref{mass}) is obvious in view of (\ref{01}) and (\ref{02}).
\qed
%
%
%
%
%
%
%
%
%
%
%
The following observation on boundedness enforced by suitably small data generalizes knowledge on similar properties
in related Keller-Segel type systems (\cite{cao_small}), and will be of
importance in our derivation both of Theorem \ref{theo41} and of Theorem \ref{theo43}.
For simplicity in presentation, we confine ourselves here to an argument based on uniform smallness of the initial data,
but we at least note that, in fact, at the cost of additional technical expense the norm appearing in \eqref{42.1} could be replaced by that in $L^\frac{n}{2}(\Omega)$.
\begin{lem}\label{lem42}
  Let $n\ge 3$ and $\Omega\subset\R^n$ be a bounded domain with smooth boundary, and let $k\ge 0$.
  Then there exists $\delta>0$ with the property that whenever $u_0\in C^0(\bom)$ is nonnegative with
  \be{42.1}
	\|u_0\|_{L^\infty(\Omega)} < \delta,
  \ee
  the solution $(u,v)$ of \eqref{01}-\eqref{02} is global and satisfies \eqref{A1} with some $C>0$.
\end{lem}
\proof
  In view of a known result from parabolic regularity theory (\cite[Theorem VI.40]{Lieberman}), it is sufficient to find $\delta>0$ such that whenever \eqref{42.1} holds, we have
  \be{42.11}
	\sup_{t\in (0,\tm)} \|\nabla v(\cdot,t)\|_{L^\infty(\Om)} <\infty.
  \ee
  To achieve this, we fix any $p>n$ and then invoke standard elliptic regularity (\cite[Thm. 19.1]{friedman}) to obtain $c_1>0$ such that
  \be{42.2}
	\|\nabla \varphi\|_{L^\infty(\Om)}^2 \le c_1 \|\-\Delta\varphi+k\varphi\|_{L^p(\Om)}^2
	\qquad \mbox{for all } \varphi\in W^{2,p}(\Om)\cap W_0^{1,p}(\Om),
  \ee
  while according to a Poincar\'e inequality (\cite[Cor. 9.1.4]{jost}, \cite[Lemma 9.1]{LW_NoDEA}) we can pick $c_2>0$ fulfilling
  \be{42.3}
	\io \varphi^2 \le c_2 \io |\nabla\varphi|^2
	\qquad \mbox{for all $\varphi\in W^{1,2}(\Omega)$ such that } \big| \{\varphi=0\}\big| \ge \frac{|\Omega|}{2}.
  \ee
  We then abbreviate
  \bas
	c_3:=\frac{2(p-1)}{pc_2},
	\quad
	c_4:=2^{\frac{2}{p}+p+1} p(p-1)c_1
	\quad \mbox{and} \quad
	c_5:=2^{\frac{2}{p}-1} p(p-1) c_1 \cdot (2^{2p} |\Om|^{1-p} )^\frac{p+2}{p},
  \eas
  and let
  \be{42.33}
	\delta:=\min \bigg\{ \Big(\frac{c_3 \oy}{2c_5}\Big)^\frac{1}{p+2} \, , \, \Big(\frac{\oy}{|\Om|}\Big)^\frac{1}{p} \bigg\}
  \ee
  with
  \bas
	\oy:=\Big(\frac{c_3}{2c_4}\Big)^\frac{p}{2},
  \eas
  observing that the first restriction in \eqref{42.33} guarantees that
  \be{42.34}
	c_3  \oy - c_4 \oy^\frac{p+2}{p} - c_5 \delta^{p+2}
	= \frac{c_3}{2} \oy \cdot \Big( 1-\frac{2c_4}{c_3} \oy^\frac{2}{p}\Big)
	+ \frac{c_3}{2} \cdot \Big( \oy - \frac{2c_5}{c_3} \delta^{p+2}\Big)
	\ge 0.
  \ee
  Now assuming $u_0\in C^0(\bom)$ to be nonnegative and such that \eqref{42.1} holds, we may use that $p>n\ge 2$, and that writing
  $a:=\frac{2}{|\Om|}\io u_0$ we thus know that $0\le \xi\mapsto (\xi-a)_+^p \in C^2([0,\infty))$, to see relying on
  \eqref{01}, Young's inequality, and \eqref{42.2} that $y(t):=\io (u(\cdot,t)-a)_+^p$, $t\in [0,\tm)$,
  belongs to $C^0([0,\tm)) \cap C^1((0,\tm))$ with
  \bea{42.4}
	y'(t) + \frac{2(p-1)}{p} \io \Big|\nabla (u-a)_+^\frac{p}{2}\Big|^2
	&=& - \frac{p(p-1)}{2} \io (u-a)_+^{p-2} |\nabla u|^2
	+ p(p-1) \io u (u-a)_+^{p-2} \nabla u\cdot\nabla v \nn\\
	&\le& \frac{p(p-1)}{2} \io u^2 (u-a)_+^{p-2} |\nabla v|^2 \nn\\
	&\le& \frac{p(p-1)c_1}{2} \|u\|_{L^p(\Om)}^2 \io u^2 (u-a)_+^{p-2} \nn\\
	&\le& \frac{p(p-1)c_1}{2} \|u\|_{L^p(\Om)}^{p+2}
	\qquad \mbox{for all } t\in (0,\tm).
  \eea
  Since \eqref{mass} ensures that $m=\io u \ge a \cdot |\{u>a\}|$ and thus $|\{u\le a\}| \ge \frac{|\Om|}{2}$ for all $t\in (0,\tm)$
  according to our choice of $a$, we may hence utilize \eqref{42.3} to estimate
  \bas
	\frac{2(p-1)}{p} \io \Big|\nabla (u-a)_+^\frac{p}{2}\Big|^2
	\ge \frac{2(p-1)}{pc_2} \io (u-a)_+^p
	= c_3 y(t)
	\qquad \mbox{for all } t\in (0,\tm),
  \eas
  whereas noting that $a\le \frac{2\delta}{|\Om|}$ by \eqref{42.1} we obtain the inequality
  \bas
	\frac{p(p-1)c_1}{2} \|u\|_{L^p(\Om)}^{p+2}
	&=& \frac{p(p-1)c_1}{2} \cdot \bigg\{ \int_{\{u\ge 2a\}} u^p + \int_{\{u<2a\}} u^p \bigg\}^\frac{p+2}{p} \\
	&\le& \frac{p(p-1)c_1}{2} \cdot \bigg\{ 2^p \int_{\{u\ge 2a\}} (u-a)^p + (2a)^p |\Om| \bigg\}^\frac{p+2}{p} \\
	&\le& \frac{p(p-1)c_1}{2} \cdot \Big\{ 2^p y(t) + 2^{2p} |\Om|^{1-p} \delta^p \Big\}^\frac{p+2}{p} \\
	&\le& 2^{\frac{2}{p}-1} p(p-1) c_1 \cdot \Big\{ (2^py(t))^\frac{p+2}{p} + (2^{2p} |\Om|^{1-p} \delta^p)^\frac{p+2}{p} \Big\} \\
	&=& c_4 y^\frac{p+2}{p}(t) + c_5 \delta^{p+2}
	\qquad \mbox{for all } t\in (0,\tm).
  \eas
  Therefore, \eqref{42.4} implies that
  \bas
	y'(t) + c_3 y(t) - c_4 y^\frac{p+2}{p}(t) - c_5\delta^{p+2} \le 0
	\qquad \mbox{for all } t\in (0,\tm),
  \eas
  so that since \eqref{42.1} along with the second requirement on $\delta$ in \eqref{42.33} guarantees that
  \bas
	y(0)= \io (u_0-a)_+^p \le \delta^p |\Om| \le \oy,
  \eas
  a comparison argument on the basis of \eqref{42.34} asserts that $y(t)\le\oy$ for all $t\in (0,\tm)$.
  As thus $\sup_{t\in (0,\tm)} \|u(\cdot,t)\|_{L^p(\Om)}$ is finite, once again relying on \eqref{42.2} we obtain \eqref{42.11}
  and conclude as intended.
\qed
\mysection{Mass bounds for steady states. Proofs of Theorem \ref{theo41} and of Corollary \ref{cor411}}
Our strategy toward proving Theorem \ref{theo41} will be based on the link between solutions to \eqref{01}-\eqref{02} and
solutions of the corresponding stationary problem
\be{stat}
 	\left\{ \begin{array}{ll}
 	\frac{\nabla u}{u}- \nabla v=0,
 	\qquad & x\in \Omega, \\[1mm]
 	\Delta v -kv+u=0,
 	\qquad & x\in \Omega, \\[1mm]
 	v=0,
 	\qquad & x\in\pO,
 	\end{array} \right.
\ee
as established through an energy-based argument in the following.
\begin{lem}\label{lem35}
  Let $n\ge 2$ and $\Omega\subset\R^n$ be a bounded domain with smooth boundary, and let $k\ge 0$ and
  $0\le u_0\in C^0(\bom)$ be such that the solution $(u,v)$ of \eqref{01}-\eqref{02} from Proposition \ref{prop33}
  is global in time and bounded in the sense that $u\in L^\infty(\Omega\times (0,\infty))$.
  Then there exist $(t_j)_{j\in\N} \subset (1,\infty)$ and functions $u_\infty$ and $v_\infty$ from
  $C^2(\bom)$ such that $u_\infty>0$ and $v_\infty\ge 0$ in $\bom$, that $t_j\to\infty$,
  $u(\cdot,t_j) \to u_\infty$ and
  $v(\cdot,t_j) \to v_\infty$ in $C^0(\bom)$ as $j\to \infty$, and that $(u_\infty,v_\infty)$ solves \eqref{stat}
  with $\io u_\infty=\io u_0$.
\end{lem}
\proof
  Using that $u>0$ in $\bom\times (0,\infty)$ by the strong maximum principle, by means of a standard computation we obtain
  the identity
  \be{35.1}
	\F(t) + \int_1^t \D(\tau) d\tau = \F(1)
	\qquad \mbox{for all } t>1,
  \ee
  where we have set
  $\F(t):=\frac{1}{2} \io |\nabla v(\cdot,t)|^2 + \frac{k}{2} \io v^2(\cdot,t) - \io u(\cdot,t) v(\cdot,t) 
  + \io u(\cdot,t) \ln u(\cdot,t)$ and
  $\D(t):=\io |2\nabla \sqrt{u(\cdot,t)} - \sqrt{u(\cdot,t)} \nabla v(\cdot,t)|^2$ for $t>0$.
  Now since $u$ is bounded and nonnegative, it readily follows that $\inf_{t>1} \F(t)>-\infty$, by \eqref{35.1} meaning that
  $\int_1^\infty \D(\tau)d\tau$ is finite, so that we can pick $(t_j)_{j\in\N} \subset (1,\infty)$ such that
  $t_j\to\infty$ and
  \be{35.2}
	2\nabla \sqrt{u(\cdot,t_j)} - \sqrt{u(\cdot,t_j)} \nabla v(\cdot,t_j) \to 0
	\qquad \mbox{a.e.~in } \Omega
  \ee
  as $j\to\infty$.
  Once more due to the boundedness of $u$, we may next invoke elliptic regularity theory (\cite{GT})
  to see that also
  $\nabla v$ is bounded in $\Omega\times (0,\infty)$, and that thus we may employ a standard result on H\"older continuity
  in parabolic equations under no-flux boundary conditions (\cite{porzio-vespri}) to obtain $\theta_1\in (0,1)$
  such that $(u(\cdot,t))_{t>1}$ is bounded in $C^{\theta_1}(\bom)$.
  Again by elliptic estimates, this entails boundedness of $(v(\cdot,t))_{t>1}$ even in $C^{2+\theta_1}(\bom)$,
  whence the Arzel\`a--Ascoli theorem provides a subsequence of $(t_j)_{j\in\N}$, for convenience again denoted by
  $(t_j)_{j\in\N}$, such that $u(\cdot,t_j) \to u_\infty$ in $C^{\theta_2}(\bom)$ and $v(\cdot,t_j) \to v_\infty$
  in $C^2(\bom)$ as $j\to\infty$ with $\theta_2:=\frac{\theta_1}{2}$ and some nonnegative limit functions
  $u_\infty\in C^{\theta_1}(\bom)$ and $v_\infty \in C^2(\bom)$ for which using \eqref{01} and \eqref{02} we can easily verify that
  $-\Delta v_\infty + k v_\infty=u_\infty$ in $\Omega$ with $v_\infty=0$ on $\pO$, and that $\io u_\infty=\io u_0$.
  Moreover, along with \eqref{35.2} this entails that as $j\to\infty$ we have
  \bas
	2\nabla \sqrt{u(\cdot,t_j)} \to \sqrt{u_\infty}\, \nabla v_\infty
	\qquad \mbox{in } C^{\theta_3}(\bom)
  \eas
  for some $\theta_3\in (0,1)$. Therefore, $\sqrt{u(\cdot,t_j)} \to \sqrt{u_\infty}$ in $C^{1+\theta_3}(\bom)$
  as $j\to\infty$ and $2\nabla \sqrt{u_\infty} \equiv \sqrt{u_\infty} \nabla v_\infty$ in $\Omega$, which in particular
  means that if we pick $x_0\in \bom$ such that $u_\infty(x_0)=\|u_\infty\|_{L^\infty(\Omega)} 
  \ge \frac{1}{|\Omega|} \io u_0>0$,
  then in the connected component $\C$ of $\{x\in\bom \ | \ u_\infty(x)>0\}$ containing $x_0$ we have
  $\nabla (\ln u_\infty-v_\infty) \equiv 0$ and hence can find $c_1>0$ such that
  $\ln u_\infty \equiv v_\infty + c_1$ in $\C$.
  As $\ln \xi\to -\infty$ as $\xi\searrow 0$, however, this ensures that actually $\C=\bom$ and that thus
  $u_\infty\equiv e^{v_\infty+c_1}$ is positive in $\bom$ and belongs to $C^2(\bom)$, and that also the first equation
  in \eqref{stat} holds throughout $\bom$.
\qed
Now a crucial observation, generalizing and quantitatively sharpening a statement from \cite{BFR} concentrating on radial
solutions in a disk, rules out large-mass steady states in strictly star-shaped two- or higher-dimensional domains:
\begin{lem}\label{lem1}
  Let $n\ge 2$ and $\Omega\subset\R^n$ be a bounded domain with smooth boundary such that
  \be{om}
	\gamma:=\min_{x\in\pO} x\cdot\nu(x)>0,
  \ee
  and suppose that $k\ge 0$.
  Then whenever $u\in C^1(\bom)\cap C^2(\Omega)$ and $v\in C^0(\bom)\cap C^2(\Omega)$ are such that
  $u>0$ and $v\ge 0$ in $\bom$ and that
  $(u,v)$ solves \eqref{stat}, we necessarily have
  \be{1.1}
	\io u \le \frac{2n|\pO|}{\gamma} + 2k|\Omega|.
  \ee
%
\end{lem}
\proof
  We firstly integrate the second equation in \eqref{stat} to see that
  \be{1.2}
	\io u = k \io v - \int_{\pO} \frac{\partial v}{\partial\nu},
  \ee
  and in order to estimate both summands on the right-hand side herein appropriately, we next use $x\cdot\nabla v$ 
  as a test function for the second equation in \eqref{stat} to find the identity
  \be{1.3}
	\io \Delta v (x\cdot\nabla v) 
	- k \io v (x\cdot\nabla v)
	= - \io u(x\cdot\nabla v).
  \ee
  Here following a well-known observation (\cite{quittner_souplet}), twice integrating by parts and using our definition of 
  $\gamma$ we obtain that
  \bea{1.4}
	\io \Delta v(x\cdot\nabla v)
	&=& - \io |\nabla v|^2 
	- \frac{1}{2} \io x\cdot\nabla|\nabla v|^2
	+ \int_{\pO} \frac{\partial v}{\partial\nu} (x\cdot\nabla v) \nn\\
	&=& \frac{n-2}{2} \io |\nabla v|^2 
	- \frac{1}{2} \int_{\pO} (x\cdot\nu)|\nabla v|^2
	+ \int_{\pO} \frac{\partial v}{\partial\nu} (x\cdot\nabla v) \nn\\
	&=& \frac{n-2}{2} \io |\nabla v|^2
	+ \frac{1}{2} \int_{\pO} (x\cdot\nu)|\nabla v|^2 \nn\\
	&\ge& \frac{\gamma}{2} \int_{\pO} |\nabla v|^2,
  \eea
  because $n\ge 2$, and because the properties $v|_{\pO}=0$ and $v\ge 0$ in $\Omega$ imply that on $\pO$
  we have $\nabla v=-|\nabla v|\nu$ and hence $\frac{\partial v}{\partial\nu} (x\cdot\nabla v)= (x\cdot \nu) |\nabla v|^2$.\\
  Apart from this, again due to the homogeneous Dirichlet boundary conditions satisfied by $v$. another integration by parts
  yields
  \be{1.5}
	- k \io v(x\cdot\nabla v)
	= - \frac{k}{2} \io x\cdot\nabla v^2
	= \frac{k}{2} \io (\nabla \cdot x) v^2
	= \frac{nk}{2} \io v^2,
  \ee
  and using that $u\nabla v=\nabla u$ by \eqref{stat} we infer from a final integration by parts that
  \be{1.6}
	- \io u(x\cdot\nabla v)
	= - \io x\cdot\nabla u
	= \io (\nabla\cdot x)u - \int_{\pO} (x\cdot\nu) u
	\le n \io u,
  \ee
  once more because $x\cdot\nu\ge 0$ by \eqref{om}.\abs
  Now a combination of \eqref{1.3} with \eqref{1.4}-\eqref{1.6} reveals that
  \bas
	\frac{\gamma}{2} \int_{\pO} |\nabla v|^2 + \frac{nk}{2} \io v^2 \le n\io u
  \eas
  and that hence, by Young's inequality,
  \bas
	k\io v - \int_{\pO} \frac{\partial v}{\partial\nu}
	&\le& k\io v + \int_{\pO} |\nabla v| \\
	&\le& \bigg\{ \frac{k}{4} \io v^2 + k|\Omega| \bigg\}
	+ \bigg\{ \frac{\gamma}{4n} \int_{\pO} |\nabla v|^2 
	+ \frac{n|\pO|}{\gamma}\bigg\} \\
	&\le& \frac{1}{2} \io u + k|\Omega| + \frac{n|\pO|}{\gamma}.
  \eas
  In conjunction with \eqref{1.2}, this entails \eqref{1.1}.
\qed
A combination of the latter two statements readily yields the first part of our main results:\abs
\proofc of Theorem \ref{theo41}. \quad
  Thanks to \eqref{41.1}, from Lemma \ref{lem1} when combined with Lemma \ref{lem35} and Proposition \ref{prop33}
  it immediately follows that the set in \eqref{MS} is not empty and hence
  $\MStar(\Omega,k)$ a well-defined nonnegative number which moreover satisfies the upper estimates in 
  \eqref{41.2} and \eqref{41.3}, respectively.
  The left inequality in \eqref{41.2} is obvious from Corollary B, whereas in the case $n\ge 3$, 
  positivity of $\MStar(\Omega,k)$ is an evident by-product of Lemma \ref{lem42}.
\qed
\proofc of Corollary \ref{cor411}. \quad
  Since for each $x\in\pO$ we have $\nu(x)=\frac{x}{|x|}$ and hence $x\cdot\nu(x)=R$,
  all statements are obvious from Theorem \ref{theo41}.
\qed
\mysection{A secondary critical mass phenomenon for radial solutions. Proof of Theorem \ref{theo43}}
In view of Corollary C, Corollary \ref{cor411}, and Lemma \ref{lem42}, verifying the occurrence of a genuinely secondary
critical mass phenomenon in the flavor of Theorem \ref{theo43} amounts to making sure that whenever the degradation 
parameter $k$ in \eqref{01} is positive, in any planar disk we can find global bounded radial solutions at some mass level
larger than $8\pi$.
To accomplish this, for such radial solutions $(u,v)=(u(r,t),v(r,t))$, $r\in[0,R]$, 
of \eqref{01}-\eqref{02} in $\Omega=B_R(0)\subset\R^2$ with $R>0$,
again maximally extended up to $\tm\in (0,\infty]$ in the style of Proposition \ref{prop33},
we follow the idea of \cite{JL} and \cite{biler_hilhorst_nadzieja_II} and introduce the cumulated quantities 
\be{w}
	w(s,t):=\int_0^{\sqrt{s}} \rho u(\rho,t) d\rho
	\qquad s\in [0,R^2], \ t\in [0,\tm),
\ee
and 
\be{z}
	z(s,t):=k\int_0^{\sqrt{s}} \rho v(\rho,t) d\rho
	\qquad s\in [0,R^2], \ t\in [0,\tm),
\ee
as well as
\be{w0}
	w_0(s):=\int_0^{\sqrt{s}} \rho u_0(\rho) d\rho, 
	\qquad s\in [0,R^2].
\ee
Then from the nonnegativity of $u$, and from \eqref{01} as well as \eqref{02}, it follows that $w_s\ge 0$ in 
$[0,R^2] \times [0,\tm)$ and
\be{0w}
    	\left\{ \begin{array}{ll}
	w_t = 4 s w_{ss} + 2ww_s - 2zw_s,
	\qquad & s\in (0,R^2), \ t\in (0,\tm), \\[1mm]
	w(0,t)=0, \quad w(R^2,t)=\frac{1}{2\pi} \cdot \io u_0,
	\qquad & t\in (0,\tm), \\[1mm]
	w(s,0)=w_0(s),
	& s\in (0,R^2),
 	\end{array} \right.
\ee
and the core of our strategy will consist in appropriately making use of the rightmost absorptive contribution to the first
equation herein in order to ensure that some of these solutions remain bounded in $C^1([0,R^2])$ even though
satisfying $w|_{s=R^2}>4$.
This will be achieved by means of a parabolic comparison with stationary supersolutions, to be 
constructed in Lemma \ref{lem25}, on the basis of a pointwise lower estimate for the function $z$	
which plays a central role
in this additional dissipative part, but which through \eqref{01}-\eqref{02} and \eqref{w} is linked to $w$ in a nonlocal manner.\abs
As a first step toward adequately coping with this, to be completed in Lemma \ref{lem24},
let us invoke a comparison argument to derive a fairly rough but useful lower bound for $w$.
\begin{lem}\label{lem21}
  Let $R>0$ and $\Omega = B_R(0)\subset\R^2$, let $k>0$, 
  and suppose that $u_0\in C^0_{rad}(\bom)$ is nonnegative and such that $w_0$ as in \eqref{w0} satisfies 
  \be{21.1}
	w_0(s) \ge \delta s^\beta
	\qquad \mbox{for all } s\in (0,R^2)
  \ee
  with some $\delta>0$ and some
  \be{21.2}
	\beta\ge 1 + \frac{1}{4\pi} \cdot \io u_0.
  \ee
  Then
  \be{21.3}
	w(s,t) \ge \delta s^\beta
	\qquad \mbox{for all $s\in (0,R^2)$ and } t\in (0,\tm).
  \ee
\end{lem}
\proof
  We abbreviate $m:=\io u_0$ and first observe that since 
  \bas
	k \io v=\io u + \int_{\pO} \frac{\partial v}{\partial\nu} \le m
	\qquad \mbox{for all } t\in (0,\tm)
  \eas
  according to the second equation in \eqref{01}  and \eqref{mass}, the function $z$ from \eqref{z} satisfies
  \bas
	z(s,t) \le  k \int_0^R \rho v(\rho,t) d\rho \le \frac{m}{2\pi}
	\qquad \mbox{for all $s\in (0,R^2)$ and any } t\in (0,\tm).
  \eas
  Therefore, writing
  \bas
	\uw(s,t):=\delta s^\beta,
	\qquad s\in [0,R^2], \ t\ge 0,
  \eas
  by nonnegativity of $\uw$ and $\uw_s$ we can estimate
  \bea{21.4}
	4s\uw_{ss} + 2\uw \uw_s - 2z(s,t)\uw_s
	&\ge& 4s\uw_{ss} - \frac{m}{\pi} \uw_s \nn\\
	&=& 4\beta(\beta-1) \delta s^{\beta-1} - \frac{m}{\pi} \cdot \beta \delta s^{\beta-1} \nn\\[1mm]
	&\ge& 0
	\qquad \mbox{for all $s\in (0,R^2)$ and } t\in (0,\tm),
  \eea
  because \eqref{21.2} asserts that $4\beta(\beta-1) \ge \frac{m}{\pi} \beta$.
  Since \eqref{21.1} implies that $\uw(s,0) \le w_0(s)$ for all $s\in (0,R^2)$, and that necessarily also
  $\uw(R^2,t) \le w_0(R^2)=w(R^2,t)$ for all $t\in (0,\tm)$ by \eqref{mass}, noting that
  $\uw(0,t)=0$ for all $t\in (0,\tm)$ we infer from the comparison principle in Lemma \ref{lem_cp} from the appendix that due to
  \eqref{21.4} indeed $w\ge \uw$ in $(0,R^2)\times (0,\tm)$.
\qed
As a consequence, we obtain the following statement on lower control of the mass accumulated in the disk 
$B_\frac{R}{2}(0)$ throughout evolution, uniform with respect to mass levels within any fixed interval.
\begin{cor}\label{cor211}
  Let $\Omega=B_R(0)\subset\R^2$ with some $R>0$, and let $k>0, m>0$, and $M\ge m$.
  Then there exists $C>0$ such that for all nonnegative $u_0\in C^0_{rad}(\bom)$ fulfilling
  \be{211.01}
	m \le \io u_0 \le M
  \ee
  as well as
  \be{211.1}
	\Mint_{B_r(0)} u_0 \ge \Mint_{B_R(0)} u_0
	\qquad \mbox{for all } r\in (0,R),
  \ee
  the solution $(u,v)$ of \eqref{01}-\eqref{02} satisfies
  \be{211.2}
	\int_{B_\frac{R}{2}(0)} u(\cdot,t) \ge C
	\qquad \mbox{for all } t\in (0,\tm).
  \ee
\end{cor}
\proof
  In order to apply Lemma \ref{lem21} to $\beta:=1+\frac{M}{4\pi}$ and $\delta:=\frac{m}{2\pi 
R^{2\beta}}$,
  we note that when rewritten in the variables $w$, $z$ and $s$ from \eqref{w} and \eqref{w0}, \eqref{211.1} together with \eqref{211.01}
  guarantees that
  \bas
	w_0(s) \ge \frac{ms}{2\pi R^2}
	\qquad \mbox{for all } s\in (0,R^2).
  \eas
  As $\beta>1$, namely, this entails that
  \bas
	\frac{w_0(s)}{\delta s^\beta}
	&\ge& \frac{m}{2\pi \delta R^2 s^{\beta-1}} \\
	&\ge& \frac{m}{2\pi \delta R^2 \cdot (R^2)^{\beta-1}} \\
	&=& \frac{m}{2\pi\delta R^{2\beta}} \\
	&=& 1
	\qquad \mbox{for all } s\in (0,R^2),
  \eas
  whence Lemma \ref{lem21} ensures that for $w$ as in \eqref{w} we have
  \bas
	w(s,t) \ge \delta s^\beta
	\qquad \mbox{for all $s\in (0,R^2)$ and any } t\in (0,\tm).
  \eas
  As a particular consequence, this implies that
  \bas
	\int_{B_\frac{R}{2}(0)} u(\cdot,t)
	= 2\pi \cdot w\Big(\frac{R^2}{4},t\Big) 
	\ge 2\pi \cdot \delta\Big(\frac{R^2}{4}\Big)^\beta
	\qquad \mbox{for all } t\in (0,\tm)
  \eas
  and thereby proves \eqref{211.2}.
\qed
This lemma will be combined with the following well-known result on positivity of the kernel associated with the
solution operator for the Helmholtz problem solved by $v$:
\begin{lem}\label{lem23}
  Let $\Omega=B_R(0)\subset \R^2$ with some $R>0$, and for $k>0$ let $G_k$ denote Green's function of $-\Delta +k$
  under homogeneous Dirichlet boundary conditions in $\Omega$. 
  Then $G_k(x,y)\ge 0$ for all $x\in\Omega$ and $y\in\Omega \setminus \{ x\}$, and there exists $C>0$ such that
  \bas
	G_k(x,y) \ge C
	\qquad \mbox{for all } (x,y) \in \Big( B_\frac{R}{2}(0) \times B_\frac{R}{2}(0) \Big)
	\setminus \Big\{ (\tilde x,\tilde y) \in B_\frac{R}{2}(0) \times B_\frac{R}{2}(0) \ \Big| \ \tilde x=\tilde y \Big\}.
  \eas
\end{lem}
\proof
  This can be found in \cite[Section 4.9]{WieKalKri}.
\qed
In fact, by means of a corresponding integral representation the function $v$ can be estimated from below in such a way that its cumulated version satisfies a linear lower bound in the following sense:
\begin{lem}\label{lem24}
  Let $\Omega=B_R(0)\subset\R^2$ with some $R>0$, and suppose that $k>0, m>0$, and $M\ge m$.
  Then there exists $C>0$ such that whenever $u_0\in C^0_{rad}(\bom)$ is nonnegative and satisfies (\ref{211.01}) as well as 
  (\ref{211.1}), the function $z$ given by (\ref{z}) fulfils
  \be{24.1}
	z(s,t) \ge C \cdot s
	\qquad \mbox{for all $s\in (0,R^2)$ and each } t\in (0,\tm).
  \ee
\end{lem}
\proof
  According to Corollary \ref{cor211}, we can pick $c_1>0$ such that for any choice of $u_0$ with the indicated properties
  we have
  \bas
	\int_{B_\frac{R}{2}(0)} u(\cdot,t) \ge c_1
	\qquad \mbox{for all  }t\in (0,\tm).
  \eas
  Thus, if relying on Lemma \ref{lem23} we fix $c_2>0$ such that Green's function $G_k$ of $-\Delta +k$ under homogeneous
  Dirichlet conditions in $\Omega$ satisfies $G_k(x,y) \ge c_2$ whenever $x\in B_\frac{R}{2}(0)$ and 
  $y\in B_\frac{R}{2}(0) \setminus \{x\}$, due to \eqref{01}-\eqref{02} and the nonnegativity of $G_k$ and $u$ we can estimate
  \bas
	v(x,t)
	&=& \int_\Omega G_k(x,y) u(y,t) dy \\
	&\ge& \int_{B_\frac{R}{2}(0)} G_k(x,y) u(y,t) dy \\
	&\ge& c_2 \int_{B_\frac{R}{2}(0)} u(y,t) dy \\
	&\ge& c_1 c_2
	\qquad \mbox{for all } x\in B_\frac{R}{2}(0) \mbox{ and } t\in (0,\tm).
  \eas
  By definition of $z$, this entails that
  \bas
	z(s,t)
	&=& \frac{ k}{2\pi} \int_{B_{\sqrt{s}}(0)} v(x,t) dx \\
	&\ge& \frac{ k}{2\pi} \cdot c_1 c_2 \cdot |B_{\sqrt{s}}(0)| \\
	&=& \frac{c_1 c_2 k}{2} \cdot s
	\qquad \mbox{for all } s\in \Big(0,\frac{R^2}{4}\Big) \mbox{ and } t\in (0,\tm).
  \eas
  As $z(\cdot,t)$ is nondecreasing on $(0,R^2)$ thanks to the nonnegativity of $v$, this moreover entails that
  \bas
	\frac{z(s,t)}{s} \ge \frac{\frac{c_1 c_2 k}{2} \cdot \frac{R^2}{4}}{R^2} 
	= \frac{c_1 c_2 k}{8}
	\qquad \mbox{for all } s\in \Big[ \frac{R^2}{4},R^2 \Big) \mbox{ and } t\in (0,\tm),
  \eas
  and that thus \eqref{24.1} holds with $C:=\frac{c_1 c_2 k}{8}$.
\qed
The key step in our derivation of Theorem \ref{theo43} can now be found in the following essentially explicit construction
of a stationary supersolution to \eqref{0w} that corresponds to a mass level exceeding the value $8\pi$.
\begin{lem}\label{lem25}
  Let $\Omega=B_R(0)\subset\R^2$ with some $R>0$, and let $k>0$.
  Then there exist $\om=\om(R,k)>8\pi$ and a function $\ow \in W^{2,\infty}((0,R^2))$ such that
  \be{25.001}
  \ow(0)=0
  \ee
  in addition to 
  \be{25.01}
    \ow(R^2)=\frac{\om}{2\pi} 
  \ee
  and
  \be{25.1}
	\ow(s) > \frac{\om s}{2\pi R^2}
	\qquad \mbox{for all } s\in (0,R^2),
  \ee
  and such that whenever $u_0 \in C^0_{rad}(\bom)$ is a nonnegative function for which $w_0$ from \eqref{w0} satisfies
  \be{25.2}
	\frac{4s}{R^2} \le w_0(s) \le \ow(s)
	\qquad \mbox{for all } s\in (0,R^2),
  \ee
  the solution of \eqref{01}-\eqref{02} has the property that
  \be{25.3}
	w(s,t) \le \ow(s)
	\qquad \mbox{for all $s\in (0,R^2)$ and } t\in (0,\tm)
  \ee
  with $w$ as defined in \eqref{w}, so that 
  \be{25.pi}
 	\sup_{(s,t)\in (0,R^2)\times(0,\tm)} \f{w(s,t)}{s} <\infty.
  \ee
\end{lem}
\proof
  Given $R>0$ and $k>0$, upon application of Lemma \ref{lem24} to $m:=8\pi$ and $M:=10\pi$ we obtain $c_1>0$ such that
  for arbitrary nonnegative $u_0\in C^0_{rad}(\bom)$ fulfilling \eqref{211.01} and \eqref{211.1}, the function 
  $z$ in \eqref{z} satisfies
  \be{25.4}
	z(s,t) \ge c_1 s
	\qquad \mbox{for all } s\in (0,R^2) \mbox{ and } t\in (0,\tm),
  \ee
  where without loss of generality we may assume that
  \be{25.44}
	c_1 \le \frac{4}{R^2}.
  \ee
  We next use that $\ln \frac{1}{s_0}\to + \infty$ as $s_0\searrow 0$ to fix $s_0\in (0,R^2)$ sufficiently small to ensure that
  \be{25.5}
	\frac{c_1}{2} \cdot \ln \frac{R^2}{s_0} > \frac{c_1}{2} + \frac{1}{R^2},
  \ee
  noting that the latter implies that
  \be{25.6}
	s_0^2 \cdot \int_{s_0}^{R^2} \sigma^{-2} e^{\frac{c_1}{2} (\sigma-s_0)} d\sigma >s_0.
  \ee
  Indeed, using that $e^\frac{c_1}{2}\xi \ge 1+\frac{c_1}{2} \xi$ for $\xi\ge 0$ shows that
  \bas
	s_0 \cdot \int_{s_0}^{R^2} \sigma^{-2} e^{\frac{c_1}{2}(\sigma-s_0)} d\sigma
	&\ge& s_0 \cdot \int_{s_0}^{R^2} \sigma^{-2} \cdot \Big\{ 1+ \frac{c_1}{2}(\sigma-s_0)\Big\} d\sigma \\
	&=& s_0 \cdot \Big(1-\frac{c_1}{2} s_0\Big) \cdot \Big( \frac{1}{s_0}-\frac{1}{R^2}\Big)
	+ \frac{c_1 s_0}{2} \cdot \ln \frac{R^2}{s_0} \\
	&=& 1 - \Big( \frac{c_1}{2} + \frac{1}{R^2}\Big) \cdot s_0
	+ \frac{c_1}{2R^2} s_0^2
	+ \frac{c_1 s_0}{2} \cdot \ln \frac{R^2}{s_0} \\
	&>& 1 + s_0 \cdot \bigg\{ \frac{c_1}{2} \cdot \ln \frac{R^2}{s_0} - \Big(\frac{c_1}{2}+\frac{1}{R^2}\Big) \bigg\} \\
	&>& 1
  \eas
  by \eqref{25.5}. Now \eqref{25.6} enables us to pick $b>0$ small enough such that
  \bas
	s_0^2 \cdot \int_{s_0}^{R^2} \sigma^{-2} e^{\frac{c_1}{2}(\sigma-s_0)} d\sigma > s_0+b,  
  \eas
  which in turn warrants the existence of $\eps\in (0,1)$ such that still
  \be{25.7}
	s_0^\frac{4+\eps}{2} \cdot \int_{s_0}^{R^2} \sigma^{-\frac{4+\eps}{2}} e^{\frac{c_1}{2}(\sigma-s_0)} d\sigma
	> s_0+b + \frac{\eps}{4b} \cdot (s_0+b)^2.
  \ee
  Observing that
  \bas
	\varphi(\xi):=
	s_0^\frac{4+\eps}{2} \cdot \int_{s_0}^{R^2} \sigma^{-\frac{4+\eps}{2}} e^{\frac{\xi}{2}(\sigma-s_0)} d\sigma,
	\qquad \xi>0,
  \eas
  in the limit $\xi\searrow 0$ satisfies
  \bas
	\varphi(\xi) 
	&\to& s_0^\frac{4+\eps}{2} \cdot \int_{s_0}^{R^2} \sigma^{-\frac{4+\eps}{2}} d\sigma \\
	&=& \frac{2}{2+\eps} \cdot s_0^\frac{4+\eps}{2} \cdot \Big( s_0^{-\frac{2+\eps}{2}} - R^{-2-\eps}\Big) \\
	&<& \frac{2}{2+\eps} \cdot s_0 \\
	&<& s_0 + b + \frac{\eps}{4b} \cdot (s_0+b)^2
  \eas
  due to e.g.~the monotone convergence theorem, from \eqref{25.7} we infer by means of a continuity argument that
  we can finally fix $c_2\in (0,c_1]$ such that the precise equality
  \be{25.8}
	s_0^\frac{4+\eps}{2} \cdot \int_{s_0}^{R^2} \sigma^{-\frac{4+\eps}{2}} e^{\frac{c_2}{2}(\sigma-s_0)} d\sigma
	= s_0+b+\frac{\eps}{4b} \cdot (s_0+b)^2
  \ee
  holds.\abs
  Upon these choices, we now let
  \be{25.9}
	\ow(s):=\left\{ \begin{array}{ll}
	\win(s)
	\qquad & \mbox{if } s\in [0,s_0], \\[1mm]
	\wout(s)
	\qquad & \mbox{if } s\in (s_0,R^2],
	\end{array} \right.
  \ee
  where
  \be{25.10}
	\win(s):=\frac{4s}{s+b},
	\qquad s\in [0,s_0],
  \ee
  which already ensures \eqref{25.001}, and where $\wout$ denotes the solution of the initial-value problem
  \be{25.11}
	\left\{ \begin{array}{l}
	4s\partial_s^2 \wout + 2(4+\eps) \partial_s \wout - 2c_2 s \cdot \partial_s \wout =0,
	\qquad s\in (s_0,R^2), \\[1mm]
	\wout(s_0)=\win(s_0), \quad 
	\partial_s \wout(s_0) = \partial_s \win (s_0).
	\end{array} \right.
  \ee
  Then $\ow$ evidently belongs to $C^1([0,R^2])\cap C^2([0,s_0]) \cap C^2([s_0,R^2])$, 
  and hence clearly also to $W^{2,\infty}((0,R^2))$, with 
  \be{25.12}
	\ow_s(s)=\frac{4b}{(s+b)^2}
	\quad \mbox{and} \quad
	\ow_{ss}(s)= - \frac{8b}{(s+b)^3}
	\qquad \mbox{for all } s\in (0,s_0),
  \ee
  and with an explicit integration of \eqref{25.11} showing that
  \bea{25.13}
	\ow_s(s)
	&=& \ow_s(s_0) \cdot \exp \bigg\{ \int_{s_0}^s \Big(-\frac{4+\eps}{2} \cdot \frac{1}{\sigma} 
	+ \frac{c_{2}}{2} \Big) d\sigma \bigg\} \nn\\
	&=& \frac{4b}{(s_0+b)^2} \cdot \Big(\frac{s_0}{s}\Big)^\frac{4+\eps}{2} e^{\frac{c_2}{2}(s-s_0)}
	\qquad \mbox{for all } s\in (s_0,R^2]
  \eea
  as well as
  \bea{25.14}
	\ow(s)
	&=& \ow(s_0)
	+ \frac{4b}{(s_0+b)^2} \cdot \int_{s_0}^s \Big(\frac{s_0}{\sigma}\Big)^\frac{4+\eps}{2} e^{\frac{c_2}{2}(\sigma-s_0)}
		d\sigma \nn\\
	&=& \frac{4s_0}{s_0+b}
	+ \frac{4b}{(s_0+b)^2} \cdot s_0^\frac{4+\eps}{2} \cdot \int_{s_0}^s \sigma^{-\frac{4+\eps}{2}}
	e^{\frac{c_2}{2}(\sigma-s_0)} d\sigma
	\qquad \mbox{for all } s\in (s_0,R^2].
  \eea
  In particular, \eqref{25.12} and \eqref{25.14} guarantee that thanks to \eqref{25.8},
  \bea{25.144}
	\ow(s)
	&\le& \ow(R^2) \nn\\
	&=& \frac{4s_0}{s_0+b}
	+ \frac{4b}{(s_0+b)^2} \cdot s_0^\frac{4+\eps}{2} \cdot \int_{s_0}^s \sigma^{-\frac{4+\eps}{2}}
	e^{\frac{c_2}{2}(\sigma-s_0)} d\sigma \nn\\
	&=& \frac{4s_0}{s_0+b}
	+ \frac{4b}{(s_0+b)^2} \cdot \Big\{ s_0+b+\frac{\eps}{4b}\cdot (s_0+b)^2\Big\} \nn\\
	&=& 4+\eps
	\qquad \mbox{for all } s\in [0,R^2],
  \eea
  while recalling the inequality $c_2\le c_1$ and \eqref{25.44} we directly obtain from \eqref{25.11} and \eqref{25.13} that
  \bas
	2s \ow_{ss}(s)
	&=& - (4+\eps-c_2 s) \ow_s(s) \\
	&\le& - (4+\eps- c_2 R^2) \ow_s(s) \\
	&\le& - (4-c_1 R^2) \ow_s(s) \\[0mm]
	&<& 0
	\qquad \mbox{for all } s\in (s_0,R^2)
  \eas
  and that hence, by \eqref{25.12},
  \bas
	\ow_{ss}(s) < 0
	\qquad \mbox{for all } s\in (0,R^2)\setminus \{s_0\}.
  \eas
  In conjunction with \eqref{25.144}, the latter concavity property in particular implies that indeed both \eqref{25.01}
  and \eqref{25.1} hold if we let $\om:=2\pi \cdot (4+\eps)$, where we note that our restriction $\eps<1$ warrants that
  $\om \le 10\pi=M$.
  As obviously also $\om\ge 8\pi=m$, assuming henceforth that $u_0\in C^0_{rad}(\bom)$ is nonnegative and such that
  \eqref{25.2} is valid, we firstly observe that \eqref{25.4} in fact applies to the function $z$ thereupon defined through
  \eqref{z}, whence again using that $c_2\le c_1$ we may infer from \eqref{25.14}, \eqref{25.4}, and \eqref{25.11} that
  \bas
	\ow_t - 4s\ow_{ss} - 2\ow \ow_s + 2z\ow_s
	&=& - 4s\ow_{ss} - 2\ow\ow_s + 2z\ow_s \\
	& \ge& -4s\ow_{ss} - 2(4+\eps) \ow_s + 2c_2 s \ow_s \\[1mm]
	&=& 0
	\qquad \mbox{for all } s\in (s_0,R^2) \mbox{ and } t\in (0,\tm),
  \eas
  whereas, simply by nonnegativity of $z$ and $\ow_s$, \eqref{25.12} ensures that
  \bas
	\ow_t - 4s\ow_{ss} - 2\ow \ow_s + 2z\ow_s
	&\ge& -4s\ow_{ss} - 2\ow\ow_s \\[1mm]
	&=& 0
	\qquad \mbox{for all } s\in (0,s_0) \mbox{ and } t\in (0,\tm).
  \eas
  Since clearly $w(0,t)=\ow(0,t)=0$ and $w(R^2,t)=\ow(R^2,t)=4+\eps$ for all $t\in (0,\tm)$, we may employ the comparison
  principle from Lemma \ref{lem_cp} to conclude that indeed \eqref{25.3} holds.
  Finally, \eqref{25.pi} follows from \eqref{25.001} together with boundedness of $\ow_s$ and \eqref{25.3}.
\qed
In order to prepare an appropriate conclusion on boundedness of $w_s$ from this, let us add the following observation on
a linear upper bound for $z$.
\begin{lem}\label{lem26}
  Let $n=2, R>0, \Omega=B_R(0)\subset\R^2$, and $k>0$ and let $u_0\in C^0_{rad}(\bom)$ be nonnegative
  and such that $w$ taken from \eqref{w} satisfies
  \be{26.1}
	\sup_{(s,t)\in (0,R^2) \times (0,\tm)} \frac{w(s,t)}{s} <\infty.
  \ee
  Then there exists $C>0$ such that
  \be{26.2}
	z(s,t) \le Cs
	\qquad \mbox{for all $s\in (0,R^2)$ and } t\in (0,\tm),
  \ee
  where $z$ is as in \eqref{z}.
\end{lem}
\proof 
  Utilizing \eqref{26.1}, let us define $c_1>0$ such that $\f{w(s,t)}{s}\le c_1$ for 
  all $s\in(0,R^2)$ and $t\in (0,\tm)$. Then since
  \[
 	z_s(R^2,t) = v(R,t) = 0 \qquad \text{for all } t\in(0,\tm)
  \]
  due to the Dirichlet condition on $v$ in \eqref{02}, and 
  since by \eqref{01} we moreover have
  \bas
	4s z_{ss}(s,t)
	= k (z(s,t) - w(s,t))
	\ge -k w(s,t)
	\qquad \mbox{for all } s\in (0,R^2) \mbox{ and } t\in (0,\tm)
  \eas
  due to the nonnegativity of $z$, on integration we infer that
  \bas
	z_s(s,t)
	&=& 0 
	- \int_s^{R^2} z_{ss}(\sigma,t) d\sigma \\
	&\le& \frac{k}{4} \int_s^{R^2} \frac{w(\sigma,t)}{\sigma} d\sigma 
	\quad \le \f{c_1 k R^2}4 =: c_2
	\qquad \mbox{for all $t\in (0,\tm)$ and any } s\in (0,R^2).
  \eas
  After one more integration, in view of the fact that $z(0,t)=0$ for all $t\in(0,\tm)$ this shows that
  \bas
	z(s,t) \le c_2 s
	\qquad \mbox{for all $t\in (0,\tm)$ and } s\in (0,R^2)
  \eas
  and thereby readily entails \eqref{26.2}.
\qed
Now employing a Bernstein-type argument in the style of \cite[Lemma 4.1]{ct_critmass}, we can indeed turn the outcome
of Lemma \ref{lem25} into an $L^\infty$ bound for $u$ by means of the following implication.
\begin{lem}\label{lem27}
  Let $n=2, R>0, \Omega=B_R(0)\subset\R^2$, and $k>0$ and let $0\not\equiv u_0\in C^0_{rad}(\bom)$ be nonnegative
  and such that $w$ from (\ref{w}) satisfies (\ref{26.1}).
  Then there exists $C>0$ such that
  \be{27.1}
	\|u(\cdot,t)\|_{L^\infty(\Omega)} \le C
	\qquad \mbox{for all } t\in (0,\tm).\
  \ee
\end{lem}
\proof
  In accordance with \eqref{26.1} and \eqref{26.2}, we first fix $c_1>0$ and $c_2>0$ such that 
  \[
 	w(s,t)\le c_1 s		
	\quad \text{and} \quad z(s,t)\le c_2s \qquad \text{for all } (s,t)\in(0,R^2)\times(0,\tm).
  \]
  With $τ:=\min\{1,\f12\tm\}$, the continuity properties of $u$ stated in Proposition \ref{prop33} enable us to find $c_3>0$ 
  satisfying 
  \begin{equation}\label{c3}
 	w_s(s,t)=u(\sqrt{s},t) \le c_3 \qquad \text{for all } s\in[0,R^2], t\in[0,τ], 
  \end{equation}
  and positivity of $u(\cdot,τ)$ in $\bom$, as ensured by the strong maximum principle, warrants the existence of $c_4>0$ such that 
  \[
 	c_4\le \f12u(\sqrt{s},τ)=w_s(s,τ) \qquad \text{for all } s\in[0,R^2].
  \]
  If for $c_5:=\min\{\f{2c_4}{c_2},\f{1}{2\pi}\norm[\Lom1]{u_0}\}\cdot  \exp(-\f{c_2}2 R^2)$ we let $\ww(s,t):=c_5(\exp(\f{c_2}2 s)-1)$, $s\in [0,R^2]$, $t\in [τ,\tm)$, 
  then $\ww(s,τ)\le c_4s$ for $s\in[0,R^2]$, $\ww(R^2,t)\le \f1{2\pi}\norm[\Lom1]{u_0} =w(R^2,t)$ for all $t\in[τ,\tm)$, and, 
  furthermore, $\ww(s,t)\ge c_6 s$ for all $(s,t)\in[0,R^2]\times[τ,\tm)$ with $c_6:=\f{c_2c_5}2$. Since 
  \[
 	\ww_t-4s\ww_{ss}-2\ww\ww_s+2z\ww_s \le 0 -4sc_5{\left(\f{c_2}2\right)}^2e^{\f{c_2}2s} + 0 + 2c_2sc_5\f{c_2}2e^{\f{c_2}2s}=0 		\text{ in } (0,R^2)\times(τ,\tm),
  \]
  a first comparison argument thus shows that 
  \begin{equation}\label{c6}
 	w(s,t)\ge \ww(s,t)\ge c_6 s \qquad \text{for all } (s,t)\in(0,R^2)\times[τ,\tm).
  \end{equation}
  To conclude our series of selections, we note that boundedness of $w$ and non-degeneracy of \eqref{0w} 
  in $(\f{R^2}2,R^2)\times(0,\tm)$ allows us to invoke parabolic Schauder theory in 
  the form of \cite[Thm. IV.10.1]{LSU} so as to obtain $c_7>0$ fulfilling 
  \begin{equation}\label{c7}
 	w_s(R^2,t)\le c_7 \qquad \text{for all } t\in[τ,\tm).
  \end{equation}
  For $α>1$, we now let 
  \[
 	y_{α}(s,t):= s^{α}\f{w_s^2(s,t)}{w(s,t)}, \qquad (s,t)\in(0,R^2]\times[τ,\tm),
  \]
  and observe that then \eqref{c6} ensures that letting $y_{α}(0,t)=0$ for $t\in[τ,\tm)$ extends $y_\alpha$ so as to become
  continuous in all of $[0,R^2]\times [\tau,\tm)$. Moreover, from \eqref{c6} and \eqref{c7} 
  we know that $y_{α}(R^2,t)\le R^{2(α-1)}\f{c_7^2}{c_6}$ for all $t\in[τ,\tm)$, while 
  combining \eqref{c3} with \eqref{c6} warrants that
  $y_{α}(s,τ)\le R^{2(α-1)}\f{c_3^2}{c_6}$ for all $s\in(0,R^2]$. 
  In the following, we fix $T\in (\tau,\tm)$ and let $(s_0,t_0)$ be any point at which the restriction 
  of $y=y_{α}$ to $(0,R^2)\times (τ,T]$ attains its maximum. Then 	
  \begin{equation}\label{ys}
 	0 = y_s = αs_0^{α-1}\f{w_s^2}{w}+2s_0^{α}\f{w_sw_{ss}}w - s_0^{α}\f{w_s^3}{w^2} 	
	\qquad \text{at } (s_0,t_0)
  \end{equation}
  and 
  \begin{align}\label{yss}
 	0\ge y_{ss} =& α(α-1)s_0^{α-2} \f{w_s^2}w + 4αs_0^{α-1}\f{w_sw_{ss}}w-2αs_0^{α-1}\f{w_s^3}{w^2}-5s_0^{α}\f{w_s^2w_{ss}}{w^2}
		\nn\\
	&+2s_0^{α}\f{w_{ss}^2}w+2s_0^{α}\f{w_sw_{sss}}w+2s_0^{α} \f{w_s^4}{w^3}
	\qquad \text{at } (s_0,t_0)
  \end{align}
  as well as 
  \begin{align}\label{yt}
 	0\le y_t &= 2s_0^{α}\f{w_s}w w_{st} - s_0^{α}\f{w_s^2}{w^2} w_t \nn\\
	&= 2s_0^{α}\f{w_s}w (4w_{ss}+4s_0w_{sss}+2w_s^2+2ww_{ss}-2z_sw_s-2zw_{ss}) - s_0^{α}\f{w_s^2}{w^2} (4s_0w_{ss}+2ww_s-2zw_s)\nn\\
	&= 4s_0\cdot 2s_0^{α}\f{w_sw_{sss}}w + 8s_0^{α}\f{w_s}{w}w_{ss} - 4zs_0^{α}\f{w_s}{w}w_{ss}-4s_0^{α+1}\f{w_s^2}{w^2} w_{ss}
		\nn\\ 
	&\qquad\qquad +4s_0^{α}w_sw_{ss}+2s_0^{α}\f{w_s^3}w-4s_0^{α}z_s\f{w_s^2}{w} + 2zs_0^{α}\f{w_s^3}{w^2}
	\qquad \text{at } (s_0,t_0).
  \end{align}
  Here we note that, evidently, \eqref{ys} entails that
  \[
 	w_{ss} = \f{w_s}2\kl{\f{w_s}{w}-\f{α}{s_0}} 
  \qquad \text{at } (s_0,t_0),
  \]
  whereas \eqref{yss} shows that hence 
  \begin{align*}
 	2s_0^{α}\f{w_sw_{sss}}{w} \le& - s_0^{α}\f{w_s^2}{2w}\kl{\f{w_s}{w}-\f{α}{s_0}}^2 
	-2αs_0^{α-1}\f{w_s^2}w\kl{\f{w_s}{w}-\f{α}{s_0}}\\
	&+\f52 s_0^{α}\f{w_s^3}{w^2}\kl{\f{w_s}w-\f{α}{s_0}} 
	-α(α-1)s_0^{α-2}\f{w_s^2}w + 2αs_0^{α-1}\f{w_s^3}{w^2}-2s_0^{α}\f{w_s^4}{w^3}\\
	=& \kl{-\f12+\f52-2} s_0^{α}\f{w_s^4}{w^3} 
	+\kl{1-2-\f52+2} αs_0^{α-1}\f{w_s^3}{w^2} 
	+\kl{-\f{α}2+2α-(α-1)}αs_0^{α-2}\f{w_s^2}{w}\\
	=& -\f32αs_0^{α-1}\f{w_s^3}{w^2} +α\kl{1+\f{α}2}s_0^{α-2}\f{w_s^2}{w}\qquad \text{at } (s_0,t_0).
  \end{align*}
  Inserting these latter two pieces of information into \eqref{yt}, we obtain 
  \begin{align*}
 	0\le& 4s_0 \kl{-\f32αs_0^{α-1}\f{w_s^3}{w^2} +α\kl{1+\f{α}2}s_0^{α-2}\f{w_s^2}{w}}
	+ 4s_0^{α}\f{w_s^2}{w}\kl{\f{w_s}{w}-\f{α}{s_0}} - 2zs_0^{α}\f{w_s^2}{w}\kl{\f{w_s}{w}-\f{α}{s_0}}\\
	&\quad -2s^{α+1}_0\f{w_s^3}{w^2} \kl{\f{w_s}{w}-\f{α}{s_0}} +2s_0^{α}w_s^2\kl{\f{w_s}{w}-\f{α}{s_0}}
	+2 s_0^{α}\f{w_s^3}w-4s_0^{α}z_s\f{w_s^2}{w} + 2zs_0^{α}\f{w_s^3}{w^2}\\
	&= -2s_0^{α+1}\f{w_s^4}{w^3} + \f{w_s^3}{w^2}(-6αs_0^{α}+4s_0^{α}+2αs_0^{α}) + \f{w_s^3}{w}(2s_0^{α}+2s_0^{α}) \\
	&\qquad +\f{w_s^2}w (2α(2+α)s_0^{α-1}-4αs_0^{α-1}) - 2αs_0^{α-1}w_s^2+2αzs_0^{α-1}\f{w_s^2}w-4s_0^{α}z_s\f{w_s^2}w\\
	\le& -2s_0^{α+1}\f{w_s^4}{w^3} +4s_0^{α}\f{w_s^3}w + 2α^2s_0^{α-1}\f{w_s^2}w + 2αzs_0^{α-1}\f{w_s^2}w\\
	\le &-s_0^{α+1} \f{w_s^4}{w^3} + 4 s_0^{α-1} w_s^2 w + 2α^2s_0^{α-1}\f{w_s^2}w + 2αzs_0^{α-1}\f{w_s^2}w\qquad\text{ in }(s_0,t_0),
  \end{align*}
  so that finally 
  \begin{align*}
 	y &= \f1{s_0}\f{w^2}{w_s^2} \cdot s_0^{α+1}\f{w_s^4}{w^3} \\
	&\le \f1{s_0}\f{w^2}{w_s^2} \kl{4s_0^{α-1} w_s^2 w + 2α^2s_0^{α-1}\f{w_s^2}w + 2αzs_0^{α-1}\f{w_s^2}w}\\
	&=4s_0^{α-2}w^3+ 2α^2s_0^{α-2}w + 2αs_0^{α-2}wz\\
	&\le 4c_1^3s_0^{α+1}+ 2α^2c_1s_0^{α-1} + 2αc_1c_2s_0^{α} 
	\qquad \text{at } (s_0,t_0).
  \end{align*}
  This entails that 
  \[
 	y_{α}(s,t)\le \max\set{R^{2(α-1)}\f{c_7^2}{c_6},R^{2(α-1)}\f{c_3^2}{c_6},4c_1^3s_0^{α+1}+2α^2c_1s_0^{α-1} + 2αc_1c_2s_0^{α}} 
	\text{ for all } (s,t)\in [0,R^2]\times[τ,\tm),
  \]
  whence letting $α\searrow 1$ we conclude that 
  \[
 	\sup_{s\in(0,R^2),t\in(τ,\tm)} \f{s}{w(s,t)} w_s^2(s,t) \le \max\set{\f{c_7^2}{c_6},\f{c_3^2}{c_6},4c_1^3R^4+2c_1 
	+ 2c_1c_2R^2},
  \]
  so that boundedness of $w_s$ in $(0,R^2)\times[τ,\tm)$, and thus of $u$ in $\Om\times[τ,\tm)$, results from 
  \eqref{26.1}. Together with \eqref{c3}, this concludes the proof.
\qed
The second of our main results has thereby actually been achieved already:\abs
\proofc of Theorem \ref{theo43}. \quad
  The first identities in \eqref{43.1}, (\ref{43.3}) and (\ref{43.4}) have precisely been stated in 
  Corollary C already. Both inequalities in (\ref{43.1}) are obvious by definition, and 
  in view of Corollary \ref{cor411}, (\ref{43.1}) directly implies (\ref{43.2}).\abs
  Finally, the strict inequality in (\ref{43.4}) can be verified by once more employing Lemma \ref{lem42},
  whereas that in (\ref{43.3}) can be seen as follows:
  Given $R>0$ and $k>0$, we take $\om(R,k)$ from Lemma \ref{lem25} and use that $\om(R,k)>8\pi$ in choosing any
  $m>8\pi$ such that $m<\om(R,k)$. Then simply defining
  \bas
	u_0(x):= \frac{m}{\pi R^2},
	\qquad x\in\bom,
  \eas
  we see on applying Lemma \ref{lem25} in conjunction with Lemma \ref{lem27} and Proposition \ref{prop33}
  that the corresponding maximally extended solution $(u,v)$ of (\ref{01})-(\ref{02}) indeed is global in time and bounded
  in the sense that (\ref{A1}) holds. In particular, this entails that indeed we must have
  $\mStar(2,R,k)\ge m>8\pi$ for any such $R$ and $k$.
\qed
\mysection{Consequences for and numerical observations concerning steady states}
\begin{cor}\label{cor44}
  Let $\Omega=B_R(0)\subset\R^2$ with some $R>0$, and let $k>0$.  
  Then for all $m<\mStar(2,R,k)$, there exists at least one pair $(u,v)\in (C^2(\bom))^2$ of radial functions
  with $u>0$ and $v\ge 0$ in $\bom$ which satisfy $\io u=m$ and solve the stationary problem (\ref{stat}) in the classical sense.
\end{cor}
\proof
  This is an evident consequence of Theorem \ref{theo43} when combined with Lemma \ref{lem35}.
\qed

\begin{figure}[hbt]
 \includegraphics[height = 55mm,keepaspectratio]{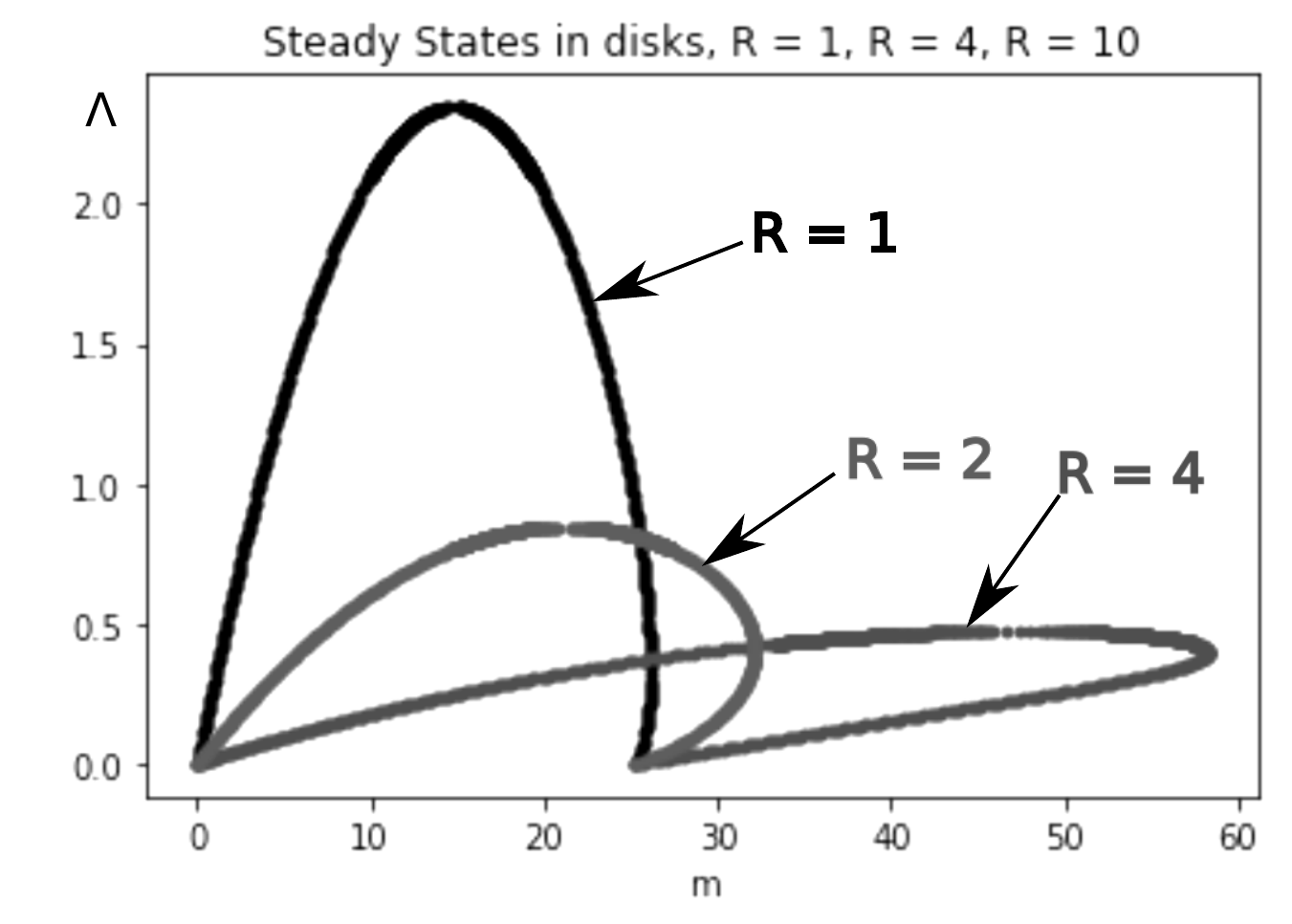} \hskip4mm %
 \includegraphics[height = 53mm,keepaspectratio]{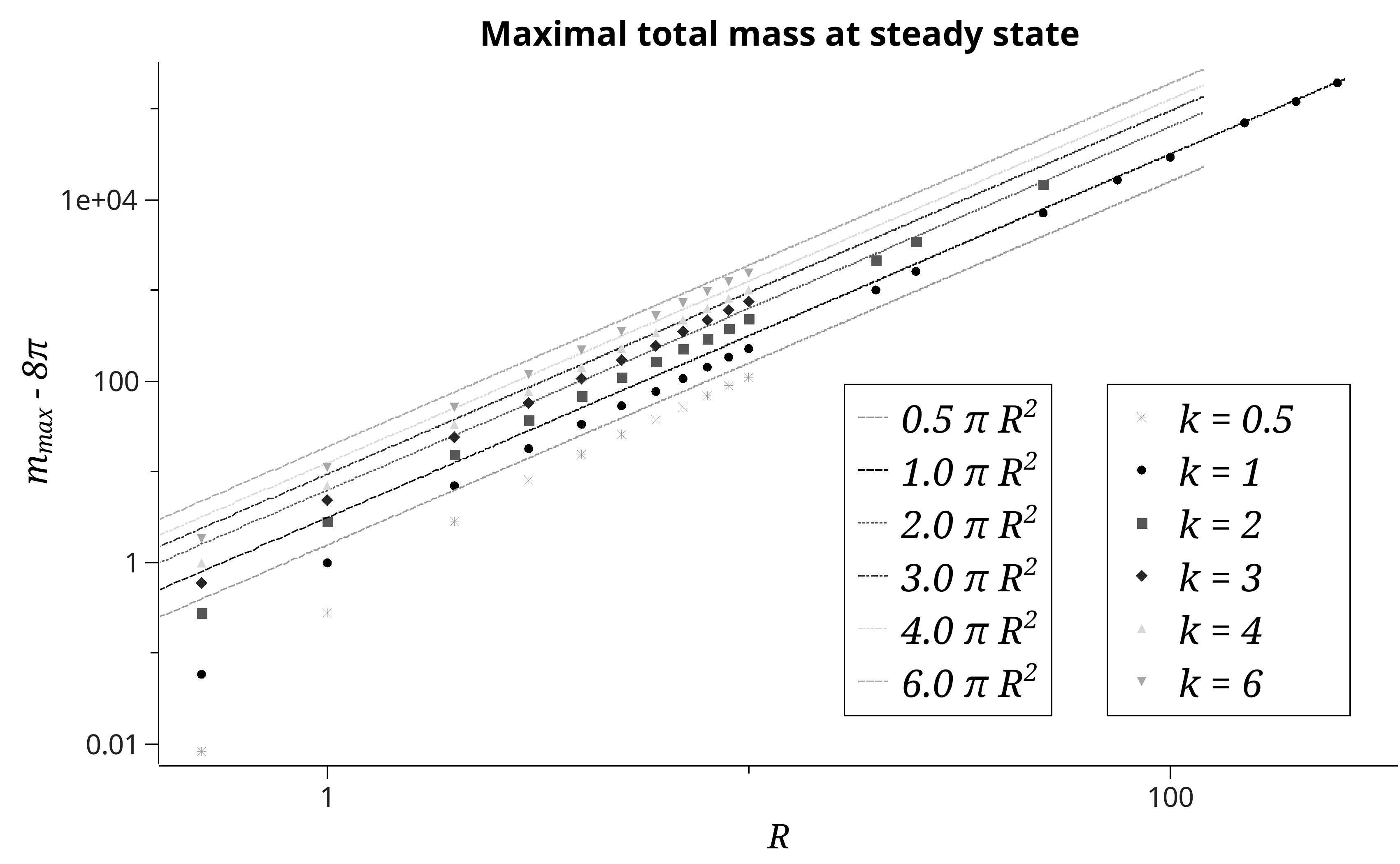}
 \caption{\small \emph{Left:} Curves of steady states as solutions of \ref{stat_scalar}; shown is the Lagrange multiplier $\Lambda$ plotted against the total mass $m=\int_\Omega u$ for $k=1$ and disks $B_R(0)\in\mathbb{R}^2$ of radii $R=1$, $R=2$, and $R=4$, respectively. Note the more pronounced tilt to the right for increasing $R$ and the common end points $(0,0)$ and $(8\pi,0)$ for all curves. \emph{Right:} Log-log plot of the maximal value of $m = \int_\Omega u$, corrected for $8\pi$, in numerically found steady state solutions in $\Omega = B_R(0) \subset\R^2$ depending on $R$ for different values of $k$. The data points are the values determined from simulation, the dashed lines correspond to the curves $m-8\pi = k\pi R^2$.}
 \label{fig:steadystates}
\end{figure}

In fact, simulations suggest the following 
\begin{conj}
 For $\Omega = B_R(0)\subset\R^2$ with $R>0$ and $k>0$, 
 \begin{enumerate}
  \item[(i) ] there is a unique steady state with $\int_\Omega u = m$ for each $m\in[0,\mstar(2,R,k)]$,
  \item[(ii) ] there are two steady states with $\int_\Omega u = m$ for each $m\in (\mstar(2,R,k),\mStar(2,R,k))$, and
  \item[(iii) ] there is a unique steady state with $\int_\Omega u = \mStar(2,R,k)$.
 \end{enumerate}
\end{conj}
As detailed in \cite{BFR}, these steady states form a continuum and can be parametrized by $\Vert u\Vert_L^\infty$. In figure \ref{fig:steadystates}, the curves of steady states in the $m$-$\Lambda$ plane are shown where $\Lambda$ is the Lagrange multiplier entering problem \eqref{stat} with $k=1$ upon integrating the first equation to $u = \Lambda \exp(v)$ and plugging this into the second equation to obtain 
\be{stat_scalar}
\left\{ \begin{array}{ll}
	-\Delta v + v= \Lambda e^v,
	\qquad & x\in \Omega, \\[1mm]
	v=0,
	\qquad & x\in\pO
	\end{array} \right.
\ee
As the curves are traced from the origin to the point $(8\pi,0)$, the norm $\Vert u\Vert_{L^\infty}$ increases, and the solution becomes more strongly concentrated near the origin. The limit point $(8\pi,0)$ would represent the singular Dirac-solution $u = 8\pi\delta_0$. The observed maximal values of $m = \int_\Omega u$ for which steady states are found, depend quadratically on the radius and hence linearly on the domain size as predicted by the upper bound on $\MStar(B_R(0),k)$ for $n=2$ from corollary \ref{cor411} and behave approximately as
\be{m_max}
m_{\max}(2,R,k) \simeq 8\pi + k\pi R^2.
\ee
\begin{figure}[hbt] \label{fig:maximizing_sols}
 \includegraphics[width=0.32\linewidth]{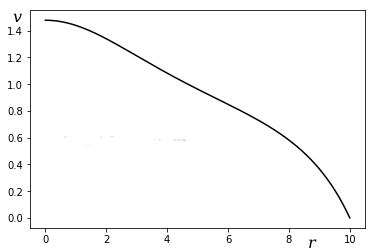} \hskip2mm %
 \includegraphics[width=0.32\linewidth]{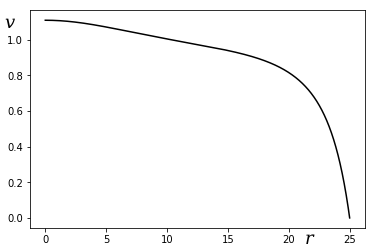} \hskip2mm %
 \includegraphics[width=0.32\linewidth]{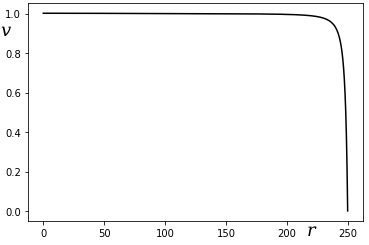} 
 \caption{Solutions $v_{\max}$ (as function of $r$) with maximal total mass $m = \Lambda\int_B e^v$ in $B_R(0)$ for $k=1$ and increasing values of $R = 10,\,25,\,250$ (left to right).}
\end{figure}
Indeed, the steady state solution maximizing the total mass for large $R$ exhibits a small peak at the origin, a wide plateau with the value $v_{plat} \approx 1$, and decreases to zero in a thin annulus given by $r\lessapprox R$. This behavior becomes obvious from the radially symmetric form 
\be{steady_state_radial}
\tilde{v}'' + \frac{\tilde{v}'}{r} - k \tilde{v} = -\Lambda \exp(\tilde{v}),~~~~~~0<r<R,~~~~~~~~~~\tilde{v}'(0)=0,~\tilde{v}(R)=0
\ee
of the steady state problem \eqref{stat_scalar}. For large $R$, the maximal value $\Lambda_c$ of $\Lambda$ allowing a solution approaches $k\,e^{-1}$, meaning that the solutions $v_\pm$ of $kv = \Lambda e^v$ are close to $1$ for $\Lambda$ close to $\Lambda_c$. Since $v_\pm$ are the values of $v$ satisfying the differential equation in \eqref{steady_state_radial} as constants, we can expect plateaus in the solution at $v\approx 1$. As we moreover observe that for large $R$ the maximal total mass is attained at $\Lambda_{\max} \lessapprox\Lambda_c$ it is not surprising that the maximal mass behaves like
\be{asympt_behavior_mmax}
m_{\max}  = 2\pi\Lambda \int_0^R r \exp(\tilde{v}(r)) dr \approx \frac{2\pi k}{e} \int_0^R r e^1 dr + \mbox{small contributions for }r\gtrapprox 0 \mbox{ and } r \approx k\pi R^2
\ee
where the small contributions of the peak near $r=0$ and the boundary layer near $r=R$ contribute with opposite signs.\abs
Figure \ref{fig:maximizing_sols} illustrates the shape of the mass maximizing solutions for different values of $R$. The plateau and lack of a pronounced peak at the origin are clearly visible for large $R=250$. 
\mysection{Discussion and biological interpretation}
Having found three distinct solvability behaviors for \eqref{01}-\eqref{02} in two dimensions, viz. global solutions for any initial conditions with $m = \int_\Omega u_0 < \Mstar(\Omega,k)$, unconditional blow up in finite time for $m>\MStar(\Omega,k)$, and the coexistence of both global and blowing up solutions for $\Mstar(\Omega,k) < m < \MStar(\Omega,k)$, we shall now briefly discuss what these results mean for the cytoskeleton of a hypothetical cell. \abs
As described in \cite{LiuPiel}, increased myosin activity -- corresponding to larger values of $m$ -- can result in the total disruption of cells. This may be interpreted as the solution to the free boundary problem associated with 
\eqref{01}-\eqref{02} (cf. \cite{BFR}) breaking down due to $\Omega$ becoming disconnected. This kind of domain blow up -- breakdown of the solution accompanied by singularities in domain shape -- has also been discussed by \cite{RPT} in one dimension where blow up in our sense -- that is, $\Vert u\Vert_{\infty} \to \infty$ in a stationary domain -- can be ruled out. Our results show that in two dimensions, the appropriate setting for a keratocyte fragment or a thin amoeboid cell on a flat substrate, classical blow-up is to be expected as well. This may be viewed as strong concentration of myosin in small regions of the cell, thereby locally disrupting the actomyosin meshwork. Clearly, upon this disruption the model will not appropriately describe the cytoskeleton anymore and would have to be replaced by another one. \abs
Acoording to this view, the regime $m<\Mstar$ will be thought of as describing a cell comfortably coming to rest on its (very sticky) substrate, and the solutions will be expected to approach the unique steady state solution with well defined distributions of myosin $u$ and the stress $v$. Increasing $m$ into the intermediate region $\Mstar < m < \MStar$ allows different fates, depending on the precise shape of the initial conditions. A cell with initial strongly concentrated myosin distribution $u_0$ will be expected to suffer disruption of its cytoskeleton while moderately concentrated $u_0$ may allow for a global solution approaching the presumably stable, weakly concentrated steady state. Further increasing $m$ beyond $\MStar$ should then lead to disruption, no matter how myosin is initially distributed inside the cell. \abs
That the difference between $\mstar(2,R,k)$ and $\mStar(2,R,k)$ increases with $R$, as suggested by figure \ref{fig:steadystates} has a physical interpretation as well. Recall that, given $k$, the cell size $R$ is measured in multiples of $\sqrt{k} \mathcal{L}$ with $\mathcal{L}$ being the viscous length of the actin gel. For small $R$, any locally generated stress will be felt throughout the cell, while for large $R$, stresses generated at one place in the cell have little impact at places far away. The stress $v$ is supposed to vanish at the boundary, and the lower branch of the two steady state solutions indicated in figure \ref{fig:steadystates} for $\mstar < m < \mStar$ comprises solutions which are monotone in $r$ but not concave down. These solutions rather feature a peak at the center of the cell, at $r=0$, where myosin is concentrated and the stress is high, a plateau at intermediate $r$ with almost constant stress and $u\approx kv$, and a region of further decreasing stress at the boundary. If the cell is large compared to the viscous length, a peak in the center can easily be established without the locally high stress being felt at the boundary, and a wider range of this type of steady states can be imagined. Recall that these steady states are expected to be unstable, and starting close to these, the solution to the time dependent problem should be expected to blow up in finite time or to relax to the supposedly stable steady state on the upper branch. \abs
It should be noted that the above discussion refers to an immobilized cell that cannot undergo shape changes or the bifurcation to a traveling wave solution. This switch from rest to steady motion occurs at even lower values $m < \Mstar$ in the free boundary problem, and it cannot be ruled out that traveling wave solutions survive as global solutions for $m>\MStar$. In fact, the local disruption of the actomyosin meshwork has been implicated in the very symmetry breaking initiating cell motion \cite{Yam}. Still, even higher values of $m$ may destroy this mode of motion and lead to physical disruption of the cell as indicated above \cite{LiuPiel}.

\mysection{Appendix: A comparison principle for \eqref{0w}}
Let us finally extract from \cite{bellomo_win2} the following comparison principle for problems of type \eqref{0w}, 
forming a reduced version of an actually more comprehensive statement involving more general degenerate parabolic operators.
\begin{lem}\label{lem_cp}
  Let $L>0$ and $T>0$, and suppose that
  $\uw$ and $\ow$ are two functions which belong to $C^1([0,L]\times [0,T))$ and satisfy
  \bas
	\uw_s(s,t)>0
	\quad \mbox{and} \quad
	\ow(s,t)>0
	\qquad \mbox{for all $s\in (0,L)$ and } t\in (0,T)
  \eas
  as well as
  \bas
	\uw(\cdot,t) \in W^{2,\infty}_{loc}((0,L))
	\quad \mbox{and} \quad
	\ow(\cdot,t) \in W^{2,\infty}_{loc}((0,L))
	\qquad \mbox{for all } t\in (0,T).
  \eas
  If for some $a\ge 0$ and some uniformly continuous 
  $b=b(s,t,\xi): (0,L) \times (0,T) \times [0,\infty)$,
  Lipschitz continuous with respect to $\xi\in [0,\xi_0]$ in $(0,L)\times (0,T)\times [0,\xi_0]$ for any $\xi_0>0$,
  we have
  \bas
 	\uw_t \le a s \uw_{ss} + b(s,t,\uw)\uw_s
	\quad \mbox{and} \quad
	\ow_t \ge a s \ow_{ss} + b(s,t,\ow)\ow_s
	\qquad \mbox{for all $t\in (0,T)$ and a.e.~$s\in (0,L)$,}
  \eas
  and if moreover
  \bas
	\uw(s,0)\le \ow(s,0)
	\qquad \mbox{for all } s\in (0,L)
  \eas
  as well as
  \bas
	\uw(0,t)\le \ow(0,t)
	\quad \mbox{and} \quad
	\uw(L,t)\le \ow(L,t)
	\qquad \mbox{for all } t\in (0,T),
  \eas
  then
  \bas
	\uw(s,t)\le \ow(s,t)
	\qquad \mbox{for all $s\in [0,L]$ and } t\in [0,T).
  \eas
\end{lem}
\proof
  This immediately results from \cite[Lemma 5.1]{bellomo_win2}.
\qed
\vspace*{5mm}
{\bf Acknowledgement.} \quad
The third author acknowledges support of the {\em Deutsche Forschungsgemeinschaft} within the 
project {\em Emergence of structures and advantages in cross-diffusion systems}, project number
411007140.\abs
The work presented here was initiated at the Mini Workshop {\it PDE models of motility and invasion in active biosystems} at the Mathematical Research Institute Oberwolfach in 2017.
{
\footnotesize 
  \setlength{\parskip}{0pt}
  \setlength{\itemsep}{0pt plus 0.2ex}

\def\cprime{$'$}

} 
\end{document}